\documentclass[reqno]{amsart}
\usepackage{graphicx}
\usepackage{color}
\usepackage{tikz}
\usepackage{khotorsion}

\usetikzlibrary{shapes.geometric}
\let\TSp\thinspace
\def\DSp{\thinspace\thinspace}
\def\lbr{\raise 1pt\hbox{[}}
\def\rbr{\raise 1pt\hbox{]}}

\def\gobble#1{}
\def\hline{\multispan\numcolumns\hrulefill\cr}
\def\torframe#1{\vtop{\vbox{\hrule\hbox{\vrule\strut #1\vrule}}\hrule}}
\def\dblhline{\hline height 0.16667em\gobble&\emptyline\cr\hline}

\def\inline#1{\leaders\hrule\hskip 0.33em plus 1fill%
  \hbox{\vbox to 0pt{\vss\hbox{\TSp #1\TSp}\vss}}%
  \leaders\hrule\hskip 0.33em plus 1fill\vrule
}

\def\knightc#1#2{%
  \lower 0.25em\hbox{\tikz[remember picture]
    \node[circle,draw,inner sep=0.1em] (#2) {#1};}%
}
\def\pawnc#1#2{%
  \lower 0.25em\hbox{\tikz[remember picture]
    \node[circle,draw,fill=black!20,inner sep=0.1em] (#2) {#1};}%
}
\def\knights#1#2{%
  \lower 0.25em\hbox{\tikz[remember picture]
    \node[regular polygon,regular polygon sides=4,draw,fill=black!20,inner sep=-0.1em] (#2) {#1};}%
}

\newbox\tablebox

\IfFileExists{\jobname.aux}{\label{KhoHo}}{\let\knights\knightc}%

\def\KhoHo{{\tt KhoHo}\space}

\begin{document}
\title{Torsion of the Khovanov homology}
\author[A.~Shumakovitch]{Alexander N. Shumakovitch}
\address{Department of Mathematics, The George Washington University,
Phillips Hall, 801\ 22nd St. NW, Suite \#739, Washington, DC 20052, U.S.A.}
\email{Shurik@gwu.edu}
\thanks{The author's initial work on this paper was partially supported
by the Swiss National Science Foundation in 2001--2003}
\begin{abstract}
Khovanov homology is a recently introduced invariant of oriented links in
$\R^3$. It categorifies the Jones polynomial in the sense that the (graded)
Euler characteristic of the Khovanov homology is a version of the Jones
polynomial for links. In this paper we study torsion of the Khovanov homology.
Based on our calculations, we formulate several conjectures about the torsion
and prove weaker versions of the first two of them. In particular, we prove
that all non-split alternating links have their integer Khovanov homology
almost determined by the Jones polynomial and signature. The only remaining
indeterminacy is that one cannot distinguish between $\Z_{2^k}$ factors in
the canonical decomposition of the Khovanov homology groups for different
values of $k$.
\end{abstract}
\maketitle

\vskip -0.65cm\null
\section{Introduction}
Let $L$ be an oriented link in the Euclidean space $\R^3$ represented by a
planar diagram~$D$. In a seminal paper, Mikhail Khovanov~\cite{Kh-Jones}
assigned to $D$ a family of abelian groups $\CalH^{i,j}(L)$, whose isomorphism
classes depend on the isotopy class of $L$ only. These groups are defined as
homology groups of an appropriate (graded) chain complex $\CalC(D)$ with
integer coefficients.
The main property of the Khovanov homology is that it {\em categorifies} the
Jones polynomial. More specifically, let $J_L(q)$ be a version of the Jones
polynomial of $L$ that satisfies the following identities (called the {\em
Jones skein relation} and {\em normalization}):
\begin{equation}\label{eq:K-jones-skein}
-q^{-2}J_{\includegraphics[scale=0.45]{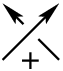}}(q)
+q^2J_{\includegraphics[scale=0.45]{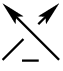}}(q)
=(q-1/q)J_{\includegraphics[scale=0.45]{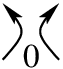}}(q);
\qquad
J_{\includegraphics[scale=0.45]{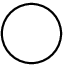}}(q)=q+1/q.
\end{equation}
The skein relation should be understood as relating the Jones polynomials of
three links whose planar diagrams are identical everywhere except in a small
disk, where they are different as depicted in~\eqref{eq:K-jones-skein}. The
normalization fixes the value of the Jones polynomial on the trivial knot.
$J_L(q)$ is a Laurent polynomial in $q$ for every link $L$ and is
completely determined by its skein relation and normalization.

The gist of the categorification is that the (graded) Euler characteristic of
the Khovanov chain complex equals $J_L(q)$:
\begin{equation}\label{eq:khpol-jones}
J_L(q)=\sum_{i,j}(-1)^iq^jh^{i,j}(L),
\end{equation}
where $h^{i,j}(L)=\rk(\CalH^{i,j}(L))$, the Betti numbers of $\CalH(L)$.
The reader is referred to~\cite{BN-first,Kh-Jones} for detailed treatment.

There are several numerical conjectures about the Khovanov homology.
We recall them briefly below.

\begin{figure}[t]
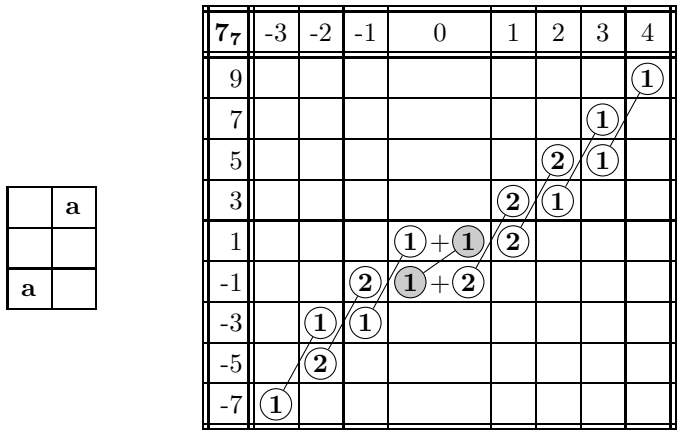

\centerline{%
\vbox{\offinterlineskip\ialign{%
\vrule height 10pt depth 5pt\enspace\hfil\bf #\hfil\enspace&
\vrule\enspace\hfil\bf #\hfil\enspace\vrule\cr
\noalign{\hrule}
&a\cr\noalign{\hrule}
&\cr\noalign{\hrule}
a&\cr\noalign{\hrule}
}\vskip 1.6cm}\qquad\qquad\input 7_7-KhHomology}%
\caption{Pattern of the ``knight move'' rule. Ranks of the Khovanov homology
of the knot $7_7$ that illustrates Conjecture~\ref{conj:barnatan-knight}.}
\label{fig:knight-move}
\end{figure}

Given a link $L$, the ranks $h^{i,j}(L)$ of its Khovanov homology can be
arranged into a table with columns and rows numbered by $i$ and $j$,
respectively (see Figure~\ref{fig:knight-move}). A pair of
entries in such a table is said to be a {\em ``knight move'' pair}, if these
entries have the same (positive) value and their $i$- and $j$-positions in the
table differ by 1 and 4, respectively. This ``knight move'' rule is depicted
in Figure~\ref{fig:knight-move}, where $a$ can be an arbitrary positive
integer.

\begin{conjec}[Bar-Natan, Garoufalidis, Khovanov~\cite{BN-first}]
\label{conj:barnatan-knight}
Let $L$ be a knot. Consider the table of Khovanov ranks $h^{i,j}(L)$ for $L$.
If one subtracts~$1$ from two adjacent entries in the column $i=0$, then the
remaining entries are arranged in ``knight move'' pairs.
\end{conjec}

\begin{example}
Figure~\ref{fig:knight-move} illustrates Conjecture~\ref{conj:barnatan-knight}
for the knot $7_7$
\footnote{Throughout this paper we use the following notation for knots:
knots with 10 crossings or less are numbered according to the Rolfsen's table
of knot~\cite{Rolfsen-book} and knots with 11 crossings or more are numbered
according to the knot table from Knotscape~\cite{Knotscape}. For example, the
knot $9_{42}$ is the knot number 42 with 9 crossings from the Rolfsen's table
and the knot $13^n_{3663}$ is a non-alternating knot number 3663 with 13
crossings from the Knotscape's one.}. The two $1$'s to subtract are shown
inside circles with a gray background and the rest of circles joined by lines
depict the ``knight move'' pairs.
\end{example}

\begin{rem}
In fact, different ``knight move'' pairs are allowed to overlap. The common
entry in this case is simply the sum of the overlapping entries from both
pairs.  For example, the knot $13^n_{3663}$, whose homology are presented in
section~\ref{app:13n-3663-knot} from the Appendix, has two overlapping pairs
$(h^{1,-1},h^{2,3})$ and $(h^{2,3},h^{3,7})$. This confusion will be cleared
up after we give a more rigorous statement of this Conjecture in
section~\ref{sec:khovanov-poly}.
\end{rem}

\begin{rem}
Conjecture~\ref{conj:barnatan-knight} was proved by Eun Soo
Lee~\cite{Lee-tor_support,Lee} for the special case of H-thin knots (see
below), in particular for all alternating knots.
\end{rem}

Let $R$ be a commutative ring with unity. In this paper, we are mainly
interested in the cases when $R=\Z$, $\Q$, or $\Z_2$.

\begin{defins}[cf. Khovanov~\cite{Kh-patterns}]
A link $L$ is said to be {\em homologically thin} over a ring $R$ or simply
{\em $R$H-thin} if its Khovanov homology groups with coefficients in $R$
are supported on two adjacent diagonals $2i-j=const$. 
A link $L$ is said to be {\em homologically slim} or simply {\em H-slim} if it
is $\Z$H-thin and all its homology groups supported on the {\em upper}
diagonal have no torsion.
A link $L$ that is not $R$H-thin is said to be {\em $R$H-thick}.
\end{defins}

\begin{prop}\label{prop:all-thins}
If a link is $\Z$H-thin, then it is $\Q$H-thin as well. Conversely, a
$\Q$H-thick link is necessarily $\Z$H-thick. If a link is H-slim, then it is
$\Z_p$H-thin for every prime $p$.
\end{prop}

\begin{figure}[t]
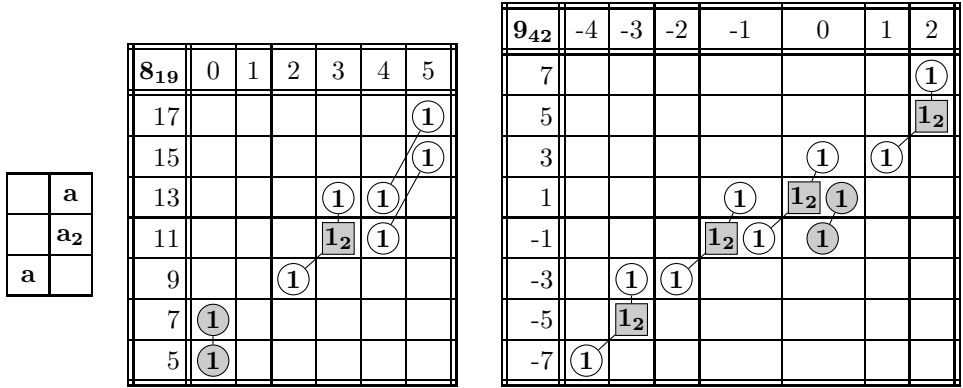

\centerline{%
\vbox{\offinterlineskip\ialign{%
\vrule height 10pt depth 5pt\enspace\hfil$\fam6#$\hfil\enspace&
\vrule\enspace\hfil$\fam6#$\hfil\enspace\vrule\cr
\noalign{\hrule}
&a\cr\noalign{\hrule}
&\kern -1em a_2\kern -1em\cr\noalign{\hrule}
a&\cr\noalign{\hrule}
}\vskip 1.2cm}\quad
\input 8_19-KhHomology\quad\input 9_42-KhHomology}%
\caption{Torsion version of the ``knight move'' rule. Khovanov homology of the
knots $8_{19}$ and $9_{42}$ that are both $\Q$H- and $\Z$H-thick but are
T-fancy and T-simple, respectively.}
\label{fig:torsion-simple}
\end{figure}

\begin{examples}
The knot $7_7$ is $\Q$H-thin since the free part of its homology is supported 
on the diagonals $2i-j=\pm1$ (see Figure~\ref{fig:knight-move}). In fact, it
is H-slim as well (see Theorem~\ref{thm:lee-thin} below). The first
$\Q$H-thick knot is $8_{19}$ (see Figure~\ref{fig:torsion-simple}).
\end{examples}

\begin{rem}
Most of the $\Z$H-thin knots are H-slim. The first prime $\Z$H-thin knot that
is not H-slim is the mirror image of $16^n_{197566}$. It is also
$\Z_2$H-thick.
\end{rem}

\begin{thm}[Lee~\cite{Lee-tor_support,Lee}]\label{thm:lee-thin}
Every oriented non-split alternating link $L$ is H-slim and the Khovanov
homology of $L$ is supported on the diagonals $2i-j=\Gs(L)\pm1$, where
$\Gs(L)$ is the signature of $L$.
\end{thm}

\begin{rem}
This theorem was originally conjectured by Bar-Natan, Garoufalidis, and
Khovanov~\cite{BN-first, Garoufalidis} in a somewhat weaker
form. Their conjecture stated that every non-split alternating link is
$\Q$H-thin and not H-slim. Lee proved a stronger version
(see~\cite{Lee-tor_support}, Corollary~4.3).
\end{rem}

Khovanov homology is very difficult to compute by hand.
Conjecture~\ref{conj:barnatan-knight} above was formulated based on extensive
computations by Dror Bar-Natan using a Mathematica software package that he
developed~\cite{BN-first}. He computed ranks of the Khovanov homology for
all prime knots with up to 11 crossings.
In~2002, the author developed \KhoHo~\cite{Sh-KhoHo} that used reductions of
the Khovanov chain complex to compute its homology faster and over $\Z$. All
the conjectures below about torsion of the Khovanov homology were formulated
based on computations done with {\tt KhoHo.} As of this writing, the integer
Khovanov homology is known for all prime knots with at most 16 crossings, all
prime links with at most 14 crossings, as well as many thousands other knots
and links.

\begin{rem}
For several years, \KhoHo was the fastest program for computing Khovanov
homology. It works efficiently for knots and links with up to 17--19
crossings. Only in Summer~2005 a significantly faster program by Bar-Natan and
Green~\cite{BN-fast,katlas-program} was written.
\end{rem}

As it turns out, torsion of the Khovanov homology has very special properties
and is at least as interesting and important as the free part of the homology.
First of all, every knot and link considered, except the unknot, the Hopf
link, their connected sums, and disjoint unions, has torsion of order~$2$. If
proved, this could lead to a (relatively) easy way to detect the unknot.

\begin{myconj}\label{conj:tors-existence}
Khovanov homology of every non-split link except the trivial knot, the Hopf
link, and their connected sums has $2$-torsion, that is, torsion elements of
order~$2$.
\end{myconj}

Torsion of orders other than $2$ appears very seldom in the Khovanov homology.
Among all 1,701,936 prime knots with at most~16 crossings, all alternating
ones have $2$-torsion only. 38 knots with~15 crossings and 129 knots with~16
crossings have $4$-torsion. One of the first such knots is the $(4,5)$-torus
knot. Its Khovanov homology is presented in section~\ref{app:4-5-torus} from
the Appendix.

\begin{rem}
The original version of this paper contained a conjecture that no link has
torsion of odd order in its Khovanov homology. This turned out to be false.
The first known counter-example is the $(5,6)$-torus knot~\cite{BN-fast} with
$3$-torsion. It has 24 crossings. Nonetheless, a large class of links is
proved to have
$2$-torsion only, see Theorem~\ref{thm:H-slim==>no-odd-torsion} below.
\end{rem}

Let $t_{p^k}^{i,j}(L)$, where $p$ is a prime number and $k\ge1$, be the
$p^k$-rank of $\CalH^{i,j}(L)$, that is, the multiplicity of $\Z_{p^k}$ in the
canonical decomposition of $\CalH^{i,j}(L)$. Let also
$T_p^{i,j}(L)=\sum_{k=1}^\infty t_{p^k}^{i,j}(L)$. The complete information
about the canonical decomposition of all the groups $\CalH^{i,j}(L)$ can be
combined into a table, where an $(i,j)$-entry contains the corresponding rank
$h^{i,j}(L)$ and a comma separated list of $t_{p^k}^{i,j}(L)$ for all relevant
$p$ and $k$ with the subscript indicating the torsion order. For example, a
hypothetical entry of $\fam6 5,3_2,1_4$ means that the corresponding group is
$\Z^5\oplus\Z_2^3\oplus\Z_4$ (see Figure~\ref{fig:torsion-simple} as well as
the Appendix). 

Similar to the case of ranks, torsion for a majority of knots and links fits
into very regular patterns that are explained below (see
Figures~\ref{fig:torsion-simple} and~\ref{fig:num-T-fancy}).
\begin{defins}\label{def:T-simple}
A link $L$ is said to be {\em weakly torsion simple} or just {\em WT-simple}
if {\bf 1)}~it satisfies Conjecture~\ref{conj:barnatan-knight}; {\bf 2)}~it
has no torsion of odd order; {\bf 3)}~for every ``knight move'' pair of value
$a$ that comprises entries $(i,j)$ and $(i+1,j+4)$, one has that
$T_2^{i+1,j+2}(L)=a$; and {\bf 4)}~all
$2^k$-ranks of $L$ fit into such patterns (see
Figure~\ref{fig:torsion-simple}, where the torsion corresponding to a ``knight
move'' pair is depicted in a gray square). $L$ is said to be {\em torsion
simple} or {\em T-simple} if it is WT-simple and has torsion of order $2$ only. It
follows that $t_2^{i+1,j+2}(L)=a$ in this case. A link $L$ that is
not WT-simple is said to be {\em T-fancy}.
\end{defins}

\begin{rem}
Strictly speaking, Conjecture~\ref{conj:barnatan-knight} is stated for knots
only. Nonetheless, it was generalized (and proved) by Lee~\cite{Lee}
to the case of $\Q$H-thin links (see Theorem~\ref{thm:lee-link-thin}).
Throughout this paper we are going to refer to a link $L$ as being T-simple or
T-fancy with understanding that this notion is assumed to be applicable, that
is, $L$ is either a knot or a $\Q$H-thin link.
\end{rem}

\begin{examples}
The knot $9_{42}$ is T-simple (see Figure~\ref{fig:torsion-simple}) and the knot
$8_{19}$ is the first T-fancy one. Both of them are $\Q$H- and $\Z$H-thick.
Theorems~\ref{thm:lee-thin} and~\ref{thm:H-slim==>WT-simple} (see below) imply
that the only T-fancy knots are non-alternating ones.
Figure~\ref{fig:num-T-fancy}
lists the number of prime T-fancy knots with at most 16 crossings.
\end{examples}

\begin{rem}
The torsion of a T-simple link is completely determined by the ranks of the
Khovanov homology. In particular, for non-split alternating links it is
completely determined by the Jones polynomial and signature~\cite{Lee}.
\end{rem}

\begin{figure}[t]
\centerline{\normalsize\def\arraystretch{1.2}%
\begin{tabular}{|l||*{9}{c|}}
\noalign{\hrule}
Number of crossings&8&9&10&11&12&13&14&15&16 \\
\noalign{\hrule\smallskip\hrule}
Number of prime&&&&&&&&&\\
\noalign{\vskip -2\smallskipamount}
\quad non-alternating knots&3&8&42&185&888&5110&27436&168030&1008906 \\
\noalign{\hrule}
Number of T-fancy knots&1&0&6&11&71&322&1736&10889&64341\\
\noalign{\hrule}
Number of T-rich knots&0&0&0&0&0&4&14&177&1413\\
\noalign{\hrule}
\end{tabular}
}
\caption{Number of prime T-fancy and T-rich knots}
\label{fig:num-T-fancy}
\end{figure}

\begin{myconj}\label{conj:H-T-simple}
Every H-slim link is T-simple. In particular, every non-split alternating link
is T-simple.
\end{myconj}

Most of the T-fancy links have their torsion ranks never greater than the
value of the corresponding ``knight move'' pair.

\begin{defin}\label{def:T-rich}
A T-fancy link $L$ is said to be {\em torsion rich} or just {\em T-rich} if it
satisfies Conjecture~\ref{conj:barnatan-knight}, has no torsion of odd order
(that is, all torsion elements have order $2^k$ for some $k$) and there is at
least one value of $(i,j)$ such that $T_2^{i+1,j+2}(L)$ is {\em greater} that
the value of the corresponding ``knight move'' pair $(h^{i,j},h^{i+1,j+4})$.
\end{defin}

\begin{example}
The first T-rich knot has 13 crossings. It is the knot $13^n_{3663}$ mentioned
above. This knot is also the first one whose homology is supported on 4
adjacent diagonals. Figure~\ref{fig:num-T-fancy} lists the number of prime
T-rich knots with at most 16 crossings.
\end{example}

In~\cite{Kh-patterns} Khovanov defined a reduced version of his
homology. The graded Euler characteristic of this reduced homology is again a
version of the Jones polynomial. More specifically, it satisfies the same skein
relation from~\eqref{eq:K-jones-skein} but is normalized to be $1$, as
opposite to $q+1/q$, on the unknot. The reduced Khovanov homology is always
supported on exactly one diagonal less than the standard one. Very few knots
have torsion in the reduced homology and all known examples of those that have
are T-rich.

\begin{myconj}\label{conj:reduced=T-rich}
A knot is T-rich if and only if its reduced Khovanov homology has torsion.
\end{myconj}

The following conjecture is a bit optimistic, but its confirmation would be
very exciting, since the torsion seems to be easier to read from the diagram
than the rest of the homology.
\begin{myconj}\label{conj:tors-ranks}
If two knots have the same torsion in their Khovanov homology, then 
they have the same ranks as well. In other words, the Khovanov homology of
a knot is completely determined by its torsion.
\end{myconj}

\begin{rem}
There are many examples of knots which have the same ranks in the Khovanov
homology but different torsion. Some of the first ones are $14^n_{9933}$ and
the mirror image of $15^n_{129763}$. The former one is T-simple, while the
latter one is T-fancy.
\end{rem}

In this paper we prove the following results.
\begin{mainthm}
\label{thm:H-slim==>no-odd-torsion}
Khovanov homology of every H-slim link has no torsion of order $p^k$ for any
odd prime $p$ and $k\ge1$.
\end{mainthm}
\begin{maincor}
Khovanov homology of every non-split alternating link has no torsion of order
$p$ for any $p$ other than a power of $2$.
\end{maincor}

We prove Theorem~\ref{thm:H-slim==>no-odd-torsion} by showing that
Conjecture~\ref{conj:barnatan-knight} holds true for every H-slim link  not
only for rational homology, but also for homology over $\Z_p$ for all odd
prime $p$. This is done using a slight modification of Lee's methods
from~\cite{Lee}.

\begin{mainthm}[cf. Conjecture~\ref{conj:H-T-simple}]
\label{thm:H-slim==>WT-simple}
Every H-slim link is WT-simple.
\end{mainthm}
\begin{maincor}
Every non-split alternating link is WT-simple. In particular, the integer
Khovanov homology of non-split alternating links is all but determined by the
Jones polynomial and signature except that one cannot distinguish between
$\Z_{2^k}$ factors in the canonical decomposition of the Khovanov homology
groups for different values of $k$.
\end{maincor}
\begin{rem}
Lee proved that the rational Khovanov homology of non-split alternating links
is completely determined by the Jones polynomial and signature, see~\cite{Lee}.
\end{rem}
\begin{maincor}[cf. Conjecture~\ref{conj:tors-existence}]
\label{cor:not-Hopf==>2-torsion}
Every alternating link except the trivial knot, the Hopf link, their connected
sums and disjoint unions has torsion of order $2$.
\end{maincor}

\begin{rem}
Marta Asaeda and Jozef Przytycki gave~\cite{Asaeda-Przytycki} an
independent proof of Corollary~\ref{cor:not-Hopf==>2-torsion}. Moreover, they
proved that an adequate link that satisfies some additional conditions has
torsion of order $2$ as well. Contrary to our approach, their proof is
constructive. They explicitly find a generator in the Khovanov chain complex
that gives rise to an appropriate torsion element. More recently, Milena
Pabiniak, Jozef Przytycki, and Radmila Sazdanovi\'c provided similar
treatments for semi-adequate
links~\cite{Jozef-Radmila+Milena-torsion,Jozef-Radmila-torsion}.
\end{rem}

Theorem~\ref{thm:H-slim==>WT-simple} is a corollary of
Theorem~\ref{thm:Z2-exact} that establishes a structure of an exact sequence
on Khovanov homology over $\Z_2$. To complete the proof of 
Conjecture~\ref{conj:H-T-simple}, one only has to show that every H-slim
link has no torsion elements of order $2^k$ for $k\ge2$.

This paper is organized as follows. Section~\ref{sec:khovanov-defs} contains
main definitions and facts about the Khovanov homology that are going to be
used in the paper. Theorems~\ref{thm:H-slim==>no-odd-torsion}
and~\ref{thm:H-slim==>WT-simple} and Corollary~\ref{cor:not-Hopf==>2-torsion}
are proved in sections~\ref{sec:khovanov-Zp},~\ref{sec:khovanov-Z2},
and~\ref{sec:2-torsion-corol-proof}, respectively. The Appendix contains
information about standard and reduced Khovanov homology of knots whose
torsion has some remarkable properties.

\subsection*{Acknowledgements}
The author is deeply grateful to Norbert A'Campo and Oleg Viro for numerous
fruitful discussions. He is also thankful to Mikhail Khovanov for many helpful
comments and for suggesting an easier proof of Theorem~\ref{thm:Z2-exact}.
Finally, the author is indebted to Jozef Przytycki for his insistence
that this paper is finished and published.

\section{Khovanov chain complex and its properties}\label{sec:khovanov-defs}

In this section we briefly recall the main ingredients of the Khovanov
homology theory. Our exposition follows the one by Viro, whose
paper~\cite{Viro-defs} is recommended for a full treatment.

\subsection{Generators and the differential of the Khovanov chain complex}
Let $D$ be a planar diagram representing an oriented link $L$. Assign a number
$\pm1$, called a {\em sign}, to every crossing of $D$ according to the rule
depicted in Figure~\ref{fig:crossing-signs}. The sum of such signs over all
the crossings is called the {\em writhe number} of $D$ and is denoted by
$w(D)$.

\begin{figure}[ht]
\captionindent 0.35\captionindent
\begin{minipage}{1.8in}
\centerline{\input{Xing_signs.pspdftex}}
\caption{Positive and negative crossings}
\label{fig:crossing-signs}
\end{minipage}
\hfill
\begin{minipage}{3.1in}
\centerline{\input{markers.pspdftex}}
\caption{Positive and negative markers and the corresponding smoothings of a
diagram.}
\label{fig:markers}
\end{minipage}
\end{figure}

At every crossing of $D$, the diagram locally divides the plane into $4$
quadrants. A choice of a pair of antipodal quadrants at a crossing can be
depicted on the diagram
with the help of a {\em marker}, which can be either {\em positive} or {\em
negative} (see Figure~\ref{fig:markers}). A collection of markers chosen at
every crossing of a diagram $D$ is called a {\em (Kauffman) state} of $D$.
There are, clearly, $2^n$ different states, where $n$ is the number of
crossings of $D$. Denote by $\Gs(s)$ the difference between the numbers of
positive and negative markers in a given state $s$.

Given a state $s$ of a diagram $D$, one can smooth $D$ at every crossing with
respect to the corresponding marker from $s$ (see Figure~\ref{fig:markers}).
The result is a family $D_s$ of disjointly embedded circles. Denote the number
of these circles by $|s|$.

Let $s$ be a state of a diagram $D$. Equip each circle from $D_s$ with either
a plus or minus sign. We call the result an {\em enhanced (Kauffman) state} of
$D$ that belongs to $s$. There are exactly $2^{|s|}$ different enhanced states
that belong to a given state $s$. Denote by $\Gt(S)$ the difference between
the numbers of positively and negatively signed circles in a given enhanced
state $S$.

With every enhanced state $S$ belonging to a state $s$ of a diagram $D$ one
can associate two numbers:
$$i(S)=\frac{w(D)-\Gs(s)}2,\qquad\qquad j(S)=-\frac{\Gs(s)+2\Gt(S)-3w(D)}2.$$
Since both $w(D)$ and $\Gs(s)$ are congruent modulo $2$ to the number of
crossings, $i(S)$ and $j(S)$ are always integer.

Fix $i,j\in\Z$. It was shown by Viro~\cite{Viro-defs} that the Khovanov chain
group $\CalC^{i,j}(D)$ is generated by all the enhanced states of $D$ with
$i(S)=i$ and $j(S)=j$. With the basis of the chain groups chosen, the Khovanov
differential $d^{i,j}:\CalC^{i,j}(D)\to\CalC^{i+1,j}(D)$ can be described by
its matrix, called the {\em incidence matrix} in this context. The elements of
the incidence matrix are called {\em incidence numbers} and are denoted by
$(S_1:S_2)$, where $S_1$ and $S_2$ are enhanced states (that is, generators)
from $\CalC^{i,j}(D)$ and $\CalC^{i+1,j}(D)$, respectively.

The incidence number $(S_1:S_2)$ is zero unless all of the following $3$
conditions are met:
\begin{enumerate}\label{def:incidence-numbers}
\renewcommand{\labelenumi}{\Roman{enumi}.}
\item the markers from $S_1$ and $S_2$ differ at one crossing of $D$ only and
at this crossing the marker from $S_1$ is positive, while the marker from
$S_2$ is negative;
\end{enumerate}
\begin{rem}
If this condition is satisfied, $D_{S_2}$ is obtained from $D_{S_1}$ by either
joining two circles into one or splitting one circle into two and, hence,
$|S_2|=|S_1|\pm1$.
\end{rem}
\begin{enumerate}
\renewcommand{\labelenumi}{\Roman{enumi}.}
\stepcounter{enumi}
\item the common circles of $D_{S_1}$ and $D_{S_2}$ have the same signs;
\item one of the following $4$ conditions are met:
\begin{enumerate}
\renewcommand{\labelenumii}{\arabic{enumii})}
\item $|S_2|=|S_1|-1$, both joining circles from $D_{S_1}$ are negative and
the resulting circle from $D_{S_2}$ is negative as well;
\item $|S_2|=|S_1|-1$, the joining circles from $D_{S_1}$ have different
signs and the resulting circle from $D_{S_2}$ is positive;
\item $|S_2|=|S_1|+1$, the splitting circle from $D_{S_1}$ is positive and
both the resulting circles from $D_{S_2}$ are positive as well;
\item $|S_2|=|S_1|+1$, the splitting circle from $D_{S_1}$ is negative and
the resulting circles from $D_{S_2}$ have different signs.
\end{enumerate}
\end{enumerate}
If all three conditions above are satisfied, then the incidence number
$(S_1:S_2)$ is defined to be equal to $(-1)^t$, where $t$ is defined as
follows. Choose some order on the crossings of $D$. Let the crossing, where
one changes the marker to get from $S_1$ to $S_2$, have number $k$ in this
order. Then $t$ is the number of negative markers in $S_1$ whose order number
is greater than $k$. As it turns out, the resulting homology does not depend
on the choice of the crossing order. More details can be found
in~\cite{BN-first, Viro-defs}.

\subsection{Reduced Khovanov homology}\label{sec:khovanov-reduced}
Let $D$ be a diagram of a link $L$. Pick a base point on $D$ that is not a
crossing. Let $\tCalC(D)$ be a subcomplex of $\CalC(D)$ generated by all the
enhanced states of $D$ that have a positive sign on the circle that the base
point belongs to. The homology $\tCalH(L)$ of this subcomplex is called the
{\em reduced Khovanov homology} of $L$. It can be shown that if $L$ is a knot,
then its reduced homology does not depend on the choice of the base point. In
general, the reduced homology of a link might depend on the component that the
base point is chosen on. See~\cite{Kh-patterns} for more details.

\begin{prop}[Khovanov~\cite{Kh-patterns}, cf.~\eqref{eq:khpol-jones}]
The graded Euler characteristic of $\tCalC(D)$ is a version of the Jones
polynomial of $L$:
\begin{equation}\label{eq:khpol-red-jones}
\tJ_L(q)=J_L(q)/(q+1/q)=\sum_{i,j}(-1)^iq^j\widetilde h^{i,j}(L).
\end{equation}
where $\widetilde h^{i,j}(L)$ are the Betti numbers of $\tCalH^{i,j}(L)$.
\end{prop}
$\tJ_L(q)$ is completely determined by the following identities
(cf.~\eqref{eq:K-jones-skein}):
\begin{equation}\label{eq:jones-skein}
-q^{-2}\tJ_{\includegraphics[scale=0.45]{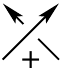}}(q)
+q^2\tJ_{\includegraphics[scale=0.45]{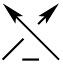}}(q)
=(q-1/q)\tJ_{\includegraphics[scale=0.45]{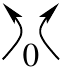}}(q);
\qquad
\tJ_{\includegraphics[scale=0.45]{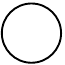}}(q)=1.
\end{equation}

\begin{prop}[Khovanov~\cite{Kh-patterns}]\label{prop:reduced-diags}
For any link $L$, its reduced Khovanov homology $\tCalH(L)$ over $\Q$ is
supported on exactly one diagonal less than the standard one.
\end{prop}
\begin{cor}\label{cor:H-thin==>alt-Jones}
If $L$ is a $\Q$H-thin link (in particular, a non-split alternating link),
then $\tCalH(L)$ is supported on exactly one diagonal. If follows that
$\tJ_L(q)$ is alternating, that is, its coefficients have alternating signs.
More precisely, if $\tJ_L(q)=\sum_{i\in\Z}c_iq^{2i+\Gg}$, where $\Gg$ is the
number of components of $L$ modulo $2$, then $\tJ_L(q)$ is alternating if and
only if $(-1)^{i-j}c_ic_j\ge0$ for all $i$ and $j$.
\end{cor}

\subsection{Khovanov polynomial and its torsion version}
\label{sec:khovanov-poly}
Let $L$ be a link and $Kh(L)(t,q)=\sum_{i,j}t^iq^jh^{i,j}(L)$ be the
Poincar\'e polynomial in variables $t$ and $q$ of its Khovanov homology. This
polynomial is called the {\em Khovanov polynomial} of $L$. Now
Conjecture~\ref{conj:barnatan-knight} can be reformulated in the following
way.

\begin{prop}[Rigorous statement of Conjecture~\ref{conj:barnatan-knight}]
Let $L$ be a knot. Then there exists a polynomial $Kh'(L)$ in $t^{\pm1}$ and
$q^{\pm1}$ with non-negative coefficients only and an even integer $s=s(L)$
such that
\begin{equation}\label{eq:knight-move}
Kh(L)=q^{s-1}\big(1+q^2+(1+tq^4)Kh'(L)\big).
\end{equation}
In other words, there exist non-negative integers $g^{i,j}(L)$ such that
\begin{equation}\label{eq:knight-move-coeffs}
h^{i,j}(L)=g^{i,j}(L)+g^{i-1,j-4}(L)+\Ge^{i,j},
\end{equation}
where $\Ge^{0,s\pm1}=1$ and $\Ge^{i,j}=0$ if $i\not=0$ or $j\not=s\pm1$.
\end{prop}

\begin{rem}
It is clear from the construction that $g^{i,j}(L)$ is the coefficient of
the term $t^iq^{j-s+1}$ in $Kh'(L)$. It must be non-zero for finitely many
values of the pair $(i,j)$ only.
\end{rem}

\begin{prop}\label{prop:H-thin-coeffs}
If $L$ is a $\Q$H-thin knot, then the polynomial $Kh'(L)$ contains powers of
$tq^2$ only. Let $Kh'(L)=\sum_{i\in Z}a_it^iq^{2i}$. In this case
$g^{i,2i+s-1}(L)=a_i$ and, hence,
\begin{alignat}{2}
h^{i,2i+s-1}(L)&=g^{i,2i+s-1}(L)+g^{i-1,2i+s-5}(L)+\Ge^{i,2i+s-1}&&
=a_i+\Gd_{i0};\\
h^{i,2i+s+1}(L)&=g^{i,2i+s+1}(L)+g^{i-1,2i+s-3}(L)+\Ge^{i,2i+s+1}&&
=a_{i-1}+\Gd_{i0},
\end{alignat}
where $\Gd_{ij}=\begin{cases} 0&\hbox{if $i\not=j$}\\ 1&\hbox{if
$i=j$}\end{cases}$ is the Kronecker delta. All other 
$h^{i,j}(L)$ are zero.
\end{prop}

Theorem~\ref{thm:lee-thin} implies that $s(L)=-\Gs(L)$ for all alternating
knots $L$, where $\Gs(L)$ is the signature of $L$.

The following theorem is a counterpart of
Conjecture~\ref{conj:barnatan-knight} for the case of $\Q$H-thin links.
\begin{thm}[Lee~\cite{Lee}]\label{thm:lee-link-thin}
Let $L$ be an $m$-component oriented $\Q$H-thin link (for example, a non-split
alternating link). Let $\ell_{k,l}$ be the linking number of the $k$-th and
$l$-th components of $L$ and let $\Gs(L)$ be the signature of $L$. Then
\begin{equation}\label{eq:lee-link-thin}
Kh(L)=q^{-\Gs(L)-1}\Biggl[(1+q^2)
\Biggl(\sum_{E\subset\{2,\dots,m\}}\!\!\!\!\!(tq^2)^
{2\sum_{\substack{\scriptscriptstyle k\in E\\\scriptscriptstyle l\not\in E}}
\ell_{k,l}}\Biggr)+(1+tq^4)Kh'(L)(tq^2)\Biggr]
\end{equation}
for some polynomial $Kh'(L)$ with non-negative coefficients.
\end{thm}

\begin{defin}
For a given link $L$, its {\em torsion Khovanov polynomial} $Kh_T$ in
variables $t^{\pm1}$ and $Q_{p^k}^{\pm1}$ is defined as
$Kh_T(L)(t,Q_2,Q_3,Q_4\,\dots)=\sum_{i,j,p,k}t^iQ_{p^k}^jt_{p^k}^{i,j}(L)$,
where $i$ and $j$ are arbitrary, $k\ge1$, and $p$ runs through all prime
numbers.
Recall that $t_{p^k}^{i,j}(L)$ is the $p^k$-rank of $\CalH^{i,j}(L)$ and
$T_p^{i,j}(L)=\sum_{k=1}^\infty t_{p^k}^{i,j}(L)$.
\end{defin}

The following proposition provides a straightforward reformulation of
Definitions~\ref{def:T-simple} and~\ref{def:T-rich} in terms of $Kh_T$.
\begin{prop}\label{prop:T-simple-coeffs}
{\bf 1.}~A link $L$ is T-simple if and only if $Kh_T(L)$ depends on the
variables $t^{\pm1}$ and $Q_2^{\pm1}$ only and $Kh_T(L)(t,q)=tq^{s+1}Kh'(L)$.
The latter equality holds true
if and only if $t_2^{i,j}(L)=g^{i-1,j-2}(L)$ for all $i$ and
$j$ (see~\ref{prop:H-thin-coeffs}).\\
{\bf 2.}~A link $L$ is WT-simple if and only if $Kh_T(L)$ depends on the
variables $t^{\pm1}$ and $Q_{2^k}^{\pm1}$ only and
$Kh_T(L)(t,q,q,\,\dots)=tq^{s+1}Kh'(L)$.
The latter equality holds true
if and only if $T_2^{i,j}(L)=g^{i-1,j-2}(L)$ for all $i$ and
$j$ (see~\ref{prop:H-thin-coeffs}).\\
{\bf 3.}~A link $L$ is T-rich if and only if $Kh_T(L)$ depends on the variables
$t^{\pm1}$ and $Q_{2^k}^{\pm1}$ only and
$tq^{s+1}Kh'(L)-Kh_T(L)(t,q,q,\,\dots)$ contains some terms with negative
coefficients.
\end{prop}

\section{Khovanov homology with $\Z_2$ coefficients}\label{sec:khovanov-Z2}
\def\Gnbar{\overline\Gn}

Denote by $\CalH_{\Z_2}^{i,j}(L)$ the Khovanov homology over $\Z_2$ (instead
of $\Z$) of an oriented link $L$ and by $h_{\Z_2}^{i,j}(L)$ the corresponding
Betti numbers. In this section we construct an acyclic differential $\Gnbar$
of bidegree $(0,2)$ on $\CalH_{\Z_2}^{i,j}(L)$ and prove
Theorem~\ref{thm:H-slim==>WT-simple}.

\subsection{Construction of the $(0,2)$-differential}
Let $D$ be a planar diagram of a link $L$ and $S$ be some enhanced state of
$D$.  Denote by $\CalN(S)$ the set of all enhanced states of $D$ that have the
same markers
and signs on all the circles as $S$ except one circle where '$+$' is
replaced with '$-$'. It is a straightforward verification that for every
enhanced state
$S'\in\CalN(S)$, one has $\Gs(S')=\Gs(S)$ and $\Gt(S')=\Gt(S)-2$. Hence,
$i(S')=i(S)$ and $j(S')=j(S)+2$.

\begin{figure}
\centerline{\scalebox{0.6}{\input{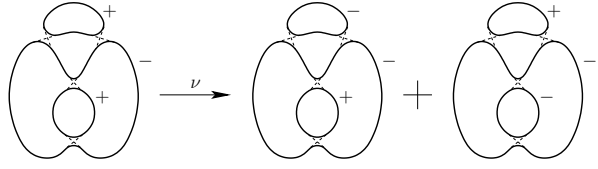}}}
\caption{Action of $\Gn$ on generators of the chain groups.}
\label{fig:nu}
\end{figure}

\begin{defin}
A differential
$\Gn^{i,j}:\CalC^{i,j}(D;\Z_2)\to\CalC^{i,j+2}(D;\Z_2)$
of bidegree $(0,2)$ on the Khovanov chain
groups with coefficients in $\Z_2$ is defined on the generators of
$\CalC^{i,j}(D;\Z_2)$ as $\Gn^{i,j}(S)=\sum_{S'\in\CalN(S)}S'$ (see
Figure~\ref{fig:nu}).
\end{defin}

To show that $\Gn$ is indeed a differential, i.e. $\Gn^2=0$, we observe that
for every enhanced states $S$ and $S''$ that have the same markers and signs
on all the circles except two circles where $S$ has '$+$', while $S''$ has
'$-$', $S''$ appears exactly twice in $\Gn(\Gn(S))$.

\begin{lem}\label{lem:phi-acyclic}
$\Gn$ is acyclic, i.e. all homology groups with respect to $\Gn$ are trivial.
\end{lem}
\begin{proof}
Let $\CalZ_n$ be the set of all length $n$ sequences of signs '$+$' and '$-$'
and let $\CalZ_n^k\subset\CalZ_n$ consists of all sequences with the difference
between the numbers of minuses and pluses being exactly $k$ (with $k\equiv
n\mod2$). $\CalZ_n$ and $\CalZ_n^k$ have $2^n$ and
$b_n^k=\dbinom{n}{\frac{n+k}2}$ elements, respectively. Denote by $G_n^k$ the
group $\Z_2^{b_n^k}$ whose factors are enumerated by the elements from
$\CalZ_n^k$.

One can construct a differential $\Gm_n^k:\CalZ_n^k\to\CalZ_n^{k+2}$ similarly
to $\Gn$: a generator of $G_n^k$ corresponding to a sequence $\Gr$ from
$\CalZ_n^k$ is mapped into the sum of the generators of $G_n^{k+2}$
corresponding to all the sequences obtained from $\Gr$ by changing exactly one
'$+$' into a '$-$'.

Let $\CalG_n$ be the complex
\begin{equation}
0\xrightarrow{}G_n^{-n}\xrightarrow{\Gm_n^{-n}}
G_n^{-n+2}\xrightarrow{\Gm_n^{-n+2}}\;\cdots\;\xrightarrow{\Gm_n^{n-4}}
G_n^{n-2}\xrightarrow{\Gm_n^{n-2}}G_n^n\xrightarrow{}0
\end{equation}
It follows from the definition of the Khovanov chain complex that the complex
\begin{equation}
\cdots\xrightarrow{\Gn^{i,j-4}}\CalC^{i,j-2}(D,\Z_2)
\xrightarrow{\Gn^{i,j-2}}\CalC^{i,j}(D,\Z_2)
\xrightarrow{\Gn^{i,j}}\CalC^{i,j+2}(D,\Z_2)
\xrightarrow{\Gn^{i,j+2}}\cdots
\end{equation}
is isomorphic to a direct sum of $\CalG_{|s|}$ with various shifts, where $s$
runs over all the Kauffman states of $D$ such that $i(s)=i$.

We claim that the complex $\CalG_n$ is acyclic. Let us prove this by induction
on $n$. The base case of $n=1$ is trivial. Denote by $\CalG_n^-$
the subcomplex of $\CalG_n$ that is obtained by choosing only those sequences
that have '$-$' it the first position. Then $\CalG_n^-$ is isomorphic to
$\CalG_{n-1}$ and, hence, is acyclic by the induction hypothesis. It follows
that $\CalG_n$ has the same homology as $\CalG_n/\CalG_n^-$. But the later is
again isomorphic to $\CalG_{n-1}$. Hence, $\CalG_n$ is acyclic as well.

This completes the proof of the lemma.
\end{proof}

\begin{lem}\label{lem:phi-commute}
$\Gn$ commutes with the Khovanov differential $d$ (over $\Z_2$).
\end{lem}

The proof is elementary and is left to the reader as an exercise.

%
%

\subsection{Patterns in $\Z_2$ homology}
Lemma~\ref{lem:phi-commute} implies that $\Gn$ can be extended to the
$(0,2)$-differential $\Gnbar$ on the Khovanov homology over $\Z_2$.
This differential is also acyclic, although this does not follow from
Lemma~\ref{lem:phi-acyclic} directly.

\begin{thm}\label{thm:Z2-exact}
$\Gnbar$ is acyclic. In particular, for every fixed $i$ the following sequence
is exact:
\begin{equation}
\cdots\xrightarrow{\Gnbar^{i,j-4}}\CalH_{\Z_2}^{i,j-2}(L)
\xrightarrow{\Gnbar^{i,j-2}}\CalH_{\Z_2}^{i,j}(L)
\xrightarrow{\Gnbar^{i,j}}\CalH_{\Z_2}^{i,j+2}(L)
\xrightarrow{\Gnbar^{i,j+2}}\cdots
\end{equation}
Consequently, $\sum_{j\in\Z}(-1)^{j}h_{\Z_2}^{i,2j+\Gg}(L)=0$ for every $i$,
where $\Gg$ is the number of components of $L$ modulo $2$.
\end{thm}
The following proof is due to Khovanov. It replaced the original one that was
too technical and unnecessarily complicated.
\begin{proof}
Choose a base point $b$ somewhere on the diagram $D$ away from the crossings.
In~\cite{Kh-patterns} Khovanov introduced another differential
$X^{i,j}:\CalC^{i,j}(D)\to\CalC^{i,j-2}(D)$ of bidegree $(0,-2)$ on the chain
groups $\CalC^{i,j}(D)$ with $\Z$ coefficients. It is defined as follows. Let
$S\in\CalC^{i,j}(D)$ be some enhanced state. If the circle of $S$ that
contains $b$ has a positive sign, then $X^{i,j}(S)=0$. Otherwise
$X^{i,j}(S)=S'\in\CalC^{i,j-2}(D)$, where $S'$ is obtained from $S$ by
changing the sign of the circle that contains $b$ from '$-$' to '$+$'. It
follows immediately from the definition that $X\circ X=0$, that is, $X$ is
indeed a differential.

It is easy to check that $X$ commutes with the Khovanov differential $d$ and,
hence, can be extended to the $(0,-2)$-differential $\overline X$ on the
Khovanov homology. We will abuse the notation slightly and denote reductions
of $X$ and $\overline X$ modulo $2$ by the same symbols.

We first claim that $\Gn\circ X+X\circ\Gn=\id$. Indeed, let $S$ be some
enhanced state of $D$. If the circle of $S$ that contains $b$ has a positive
sign, then $\Gn(X(S))=0$. Moreover, $\Gn(S)$ is a sum of enhanced states
such that all but one of them
have positive signs on their circles that contain $b$. Hence,
$X(\Gn(S))=S$. On the other hand, if the circle of $S$ containing $b$ has a
negative sign, then all the enhanced states from $\Gn(S)$ have negative signs
on their circles containing $b$ and $X(\Gn(S))$ is the sum of all the states
that are obtained from $S$ by changing the sign of the circle that contains
$b$ from '$-$' to '$+$' and changing the sign of some other circle from '$+$'
to '$-$'. Moreover, $X(S)$ has one more positive sign than $S$ and
$\Gn(X(S))$ is the sum of all the same states as $X(\Gn(S))$ plus $S$ itself.
The claim follows.

Since $\Gn$ and $X$ both commute with the differential $d$, one has that
$\Gnbar\circ\overline X+\overline X\circ\Gnbar=\id$ on the homology level 
as well. It follows that
$\Gnbar$ is acyclic. Indeed, if $\Ga\in\CalH^{i,j}(D)$ such that
$\Gnbar(\Ga)=0$, then $\Gnbar(\overline X(\Ga))=\Ga$, that is, $\Ga$ lies in
the image of $\Gnbar$.
\end{proof}

\begin{cor}
$\overline X$ is acyclic on $\CalH_{\Z_2}(L)$ as well.
\end{cor}

\begin{cor}
Since the reduced Khovanov homology over $\Z_2$
is isomorphic to the kernel of $\overline X$, it follows that 
$\CalH_{\Z_2}^{i,j}(L)\simeq\tCalH_{\Z_2}^{i,j-1}(L)\oplus\tCalH_{\Z_2}^{i,j+1}(L)$
for every $i,j$.
\end{cor}

\begin{cor}\label{cor:equal-Z2-diags}
Let $L$ be an H-slim link. Then it is $\Z_2$H-thin by~\ref{prop:all-thins},
that is, its $\Z_2$-homology are supported on the diagonals $2i-j=-s\pm1$.
Theorem~\ref{thm:Z2-exact} implies that
$h_{\Z_2}^{i,2i+s-1}(L)=h_{\Z_2}^{i,2i+s+1}(L)$ for every $i$.
\end{cor}

\begin{proof}[Proof of Theorem~\ref{thm:H-slim==>WT-simple}]
Let $L$ be an H-slim knot. It follows from
Theorem~\ref{thm:H-slim==>no-odd-torsion} that $L$ has torsion of order~$2^k$
only. It remains to show that $T_2^{i,j}(L)=g^{i-1,j-2}(L)$ for all $i$ and
$j$ (see~\ref{prop:T-simple-coeffs}). Since $T_2^{i,j}(L)=g^{i-1,j-2}(L)=0$ for
$j\not=2i+s-1$ and $g^{i-1,2i+s-3}(L)=a_{i-1}$ in notation
of~\ref{prop:H-thin-coeffs}, we only need to prove that
$T_2^{i,2i+s-1}(L)=a_{i-1}$.

Observe now that $h_{\Z_2}^{i,j}(L)=h^{i,j}(L)+T_2^{i,j}(L)+T_2^{i+1,j}(L)$.
It follows from~\ref{prop:H-thin-coeffs} that
$h_{\Z_2}^{i,2i+s-1}=a_i+\Gd_{i0}+T_2^{i,2i+s-1}$ and
$h_{\Z_2}^{i,2i+s+1}=a_{i-1}+\Gd_{i0}+T_2^{i+1,2i+s+1}$.
Corollary~\ref{cor:equal-Z2-diags} implies that
\begin{equation}
T_2^{i,2i+s-1}-a_{i-1}=T_2^{i+1,2i+s+1}-a_i
\end{equation}
for all $i$. Hence, $T_2^{i,2i+s-1}-a_{i-1}=const$, for some constant 
independent of $i$. Since the support of the Khovanov homology is finite,
there exists $i$ such that $T_2^{i,2i+s-1}=a_{i-1}=0$. It follows that the
constant must be zero.

The case of $L$ being a link can be considered similarly.
\end{proof}

\subsection{Proof of Corollary~\ref{cor:not-Hopf==>2-torsion}}
\label{sec:2-torsion-corol-proof}
Let $L$ be an alternating link with $m$ components. Then its Jones polynomial
has the form $\tJ_L(q)=\sum_i c_iq^{2i+\Gg}$, where $\Gg=m\mod 2$ (cf.
Corollary~\ref{cor:H-thin==>alt-Jones}). Define
$d(L)=|\tJ_L(\sqrt{-1})|=\sum_i|c_i|$. In fact, $d(L)=|\det(L)|$, where
$\det(L)$ is the determinant of $L$, hence the notation.

\begin{thm}[Thistlethwaite~\cite{Thistlethwaite-trees}, Theorem~1]
\label{thm:thistlethwaite}
Let $L$ be a prime non-split alternating link that admits an irreducible
alternating diagram with $n$ crossings, and let 
$\tJ_L(q)=\sum_{i=u}^v c_iq^{2i+\Gg}$ with $c_u\not=0$ and $c_v\not=0$ be its
Jones polynomial. Then $v-u=n$ and $c_ic_{i+1}\le0$ for every $u\le i<v$.
If, moreover, $L$ is not a $(2,k)$-torus link, then $c_i\not=0$ for every
$u\le i\le v$.
\end{thm}
\begin{lem}\label{lem:alt-determinant}
Let $L$ be an alternating link with $m$ components. Then $d(L)\ge2^{m-1}$.
Moreover, if $L$ is not the trivial knot, the Hopf link, their connected sums
or disjoint unions, then $d(L)>2^{m-1}$.
\end{lem}
\begin{proof}
Assume first that $L$ is non-split and prime. It is easy to deduce
from~\eqref{eq:jones-skein} (see also~\cite{Jones}) that
$\tJ_L(1)=\sum_i c_i=(-2)^{m-1}$. Hence, $d(L)\ge|\tJ_L(1)|=2^{m-1}$. It only
remains to show that if $L$ is neither the trivial knot, nor the Hopf link,
then the polynomial $\tJ_L(q)$ has both strictly positive and strictly negative
coefficients. Theorem~\ref{thm:thistlethwaite} implies this for all $L$ except
the trivial knot and $(2,k)$-torus links.

Denote a $(2,k)$-torus link by $TL_k$. It easily follows
from~\eqref{eq:jones-skein} that $d(TL_k)=k$. Indeed, if one changes a
positive crossing of $TL_k$ into a negative one, one gets $TL_{k-2}$ , and if
one smoothes such a crossing, one obtains $TL_{k-1}$. Recall now that a
$(2,k)$-torus link has at most $2$ components and that the Hopf link is
exactly the $(2,2)$-torus link.

Now the Lemma follows from the fact that the Jones polynomial of the connected
sum and disjoint union of links $L_1$ and $L_2$ is equal to
$\tJ_{L_1}(q)\tJ_{L_2}(q)$ and $(q+1/q)\tJ_{L_1}(q)\tJ_{L_2}(q)$, respectively.
\end{proof}

\begin{lem}\label{lem:alt-rank}
Let $L$ be an alternating link that is not the trivial knot, the Hopf link,
their connected sums and disjoint unions.
Then $\rank\CalH(L)>2^m$, where $m$
is the number of components of $L$ and $\rank\CalH(L)=\sum_{i,j}h^{i,j}(L)$ is
the total rank of the Khovanov homology.
\end{lem}
\begin{proof}
Consider first the case when $L$ is non-split. Theorem~\ref{thm:lee-link-thin}
implies that $\rank\CalH(L)\ge2^m$. Assume that this rank is $2^m$. In this
case $Kh'(L)$ must be $0$ and
\begin{equation}
J_L(q)=Kh(L)(-1,q)=
q^{-\Gs(L)}\Biggl[(q+1/q)\Biggl(\sum_{E\subset\{2,\dots,m\}}\!\!\!\!\!(-q^2)^
{2\sum_{\substack{\scriptscriptstyle k\in E\\\scriptscriptstyle l\not\in E}}
\ell_{k,l}}\Biggr)\Biggr].
\end{equation}
Since $\tJ_L(q)=J_L(q)/(q+1/q)$, one has
\begin{equation}
d(L)=|\tJ_L(\sqrt{-1})|=\sum_{E\subset\{2,\dots,m\}}\!\!\!\!\!1^
{2\sum_{\substack{\scriptscriptstyle k\in E\\\scriptscriptstyle l\not\in E}}
\ell_{k,l}}=2^{m-1}.
\end{equation}
This contradicts Lemma~\ref{lem:alt-determinant}. Hence, $\rank\CalH(L)>2^m$.

The general case follows from the fact that $\rank\CalH(L)$ is multiplicative
under disjoint union (see~\cite{Kh-Jones}, Corollary~12).
\end{proof}

Let us now finish the proof of Corollary~\ref{cor:not-Hopf==>2-torsion}.
Lemma~\ref{lem:alt-rank} states that $\rank\CalH(L)>2^m$ and, hence,
$Kh'(L)\not=0$ (in notation of Theorem~\ref{thm:lee-link-thin}). Since $L$ is
WT-simple by Theorem~\ref{thm:H-slim==>WT-simple}, it follows
from~\ref{prop:T-simple-coeffs} that $Kh_T(L)\not=0$ as well. Hence, $L$ has
non-trivial torsion. Since $L$ is WT-simple, some torsion elements must be of
order~$2$.
\qed

\def\vara{\mathbf{a}}
\def\varb{\mathbf{b}}

\section{Torsion of order $p$ of the Khovanov homology}\label{sec:khovanov-Zp}

This section is devoted to proving Theorem~\ref{thm:H-slim==>no-odd-torsion}.
We start by showing that the Khovanov homology over $\Z_p$ of an H-slim link
satisfy Conjecture~\ref{conj:barnatan-knight} as well.

\subsection{Khovanov homology with $\Z_p$ coefficients}
Let $L$ be an oriented link and $p$ be an odd prime number. Denote by
$\CalH_{\Z_p}^{i,j}(L)$ the Khovanov homology of $L$ over $\Z_p$ and by
$h_{\Z_p}^{i,j}(L)$ their Betti numbers. Let
$Kh_{\Z_p}(L)(t,q)=\sum_{i,j}t^iq^jh_{\Z_p}^{i,j}(L)$ be the corresponding
Poincar\'e polynomial.

\begin{thm}[cf. Theorem~\ref{thm:lee-link-thin} and~\cite{Lee},
Theorems~1.2 and~1.4]\label{thm:Zp-link-thin}
Let $L$ be an $m$-component oriented H-slim link, for example, a non-split
alternating link. Then $Kh_{\Z_p}(L)$ satisfies
identity~\eqref{eq:lee-link-thin} for the original Khovanov polynomial with
some other polynomial $Kh'_p(L)$ instead of $Kh'(L)$. If $L$ is an H-slim
knot, then this identity becomes
\begin{equation}\label{eq:knight-move-Zp}
Kh_{\Z_p}(L)=q^{-\Gs(L)-1}\bigl(1+q^2+(1+tq^4)Kh'_p(L)(tq^2)\bigr),
\end{equation}
where $\Gs(L)$ is the signature of $L$.
\end{thm}
\begin{proof}
We will show that the methods used by Lee to prove
Theorem~\ref{thm:lee-link-thin} for the Khovanov homology with $\Q$
coefficients work in our $\Z_p$ case as well if $p$ is odd prime. Only the
main steps are to be outlined and the reader is assumed to be familiar
with~\cite{Lee}.

First of all, we define the Lee differential $\GF$ of bidegree $(1,4)$ on the
Khovanov chain complex $\CalC$. The corresponding incidence numbers
$(S_1:S_2)_\GF$ of two enhanced states $S_1\in\CalC^{i,j}_{\Z_p}(D)$ and
$S_2\in\CalC^{i+1,j+4}_{\Z_p}(D)$ are defined is a similar way to the original
ones from page~\pageref{def:incidence-numbers} with the only difference being
in condition III:\\
The incidence number $(S_1:S_2)_\GF$ is zero unless all of the following $3$
conditions are met, in which case $(S_1:S_2)_\GF=\pm1$ with the sign defined
as before:
\begin{enumerate}
\renewcommand{\labelenumi}{\Roman{enumi}${}_\GF$.}
\item the markers from $S_1$ and $S_2$ differ at one crossing of $D$ only and
at this crossing the marker from $S_1$ is positive, while the marker from
$S_2$ is negative;
\item the common circles of $D_{S_1}$ and $D_{S_2}$ have the same signs;
\item one of the following two conditions are met:
\begin{enumerate}
\renewcommand{\labelenumii}{\arabic{enumii})}
\item $|S_2|=|S_1|-1$, both joining circles from $D_{S_1}$ are positive and
the resulting circle from $D_{S_2}$ is negative;
\item $|S_2|=|S_1|+1$, the splitting circle from $D_{S_1}$ is positive and
both the resulting circles from $D_{S_2}$ are negative;
\end{enumerate}
\end{enumerate}

It is easy to see that $\GF$ is indeed a differential and (anti)commutes with
the Khovanov differential $d$. Lee's proofs from~\cite{Lee} can be
applied to our version of $\GF$ to show that it also commutes with the
isomorphisms induced on $\CalH_{\Z_p}(L)$ by the Reidemeister moves. Hence,
$\GF$ gives rise to a well defined differential on $\CalH_{\Z_p}(L)$.

Consider now yet another differential $\GF+d$ on $\CalC_{\Z_p}(D)$. It can
be best described by changing the labels on the circles comprising enhanced
states from '$+$' and '$-$' to $\vara=(\hbox{'$+$'})+(\hbox{'$-$'})$ and
$\varb=(\hbox{'$+$'})-(\hbox{'$-$'})$. In this notation one has a new third
condition on the incidence numbers $(S_1:S_2)_{\GF d}$:
\begin{enumerate}
\renewcommand{\labelenumi}{\Roman{enumi}${}_{\GF d}$.}
\stepcounter{enumi}\stepcounter{enumi}
\item one of the following $4$ conditions are met:
\begin{enumerate}
\renewcommand{\labelenumii}{\arabic{enumii})}
\item $|S_2|=|S_1|-1$, both joining circles from $D_{S_1}$ and the resulting
circle from $D_{S_2}$ are marked with $\vara$, then $(S_1:S_2)_{\GF d}=\pm2$;
\item $|S_2|=|S_1|-1$, both joining circles from $D_{S_1}$ and the resulting
circle from $D_{S_2}$ are marked with $\varb$, then $(S_1:S_2)_{\GF d}=\mp2$;
\item $|S_2|=|S_1|+1$, the splitting circle from $D_{S_1}$ and both the
resulting circles from $D_{S_2}$ are marked with $\vara$, then $(S_1:S_2)_{\GF
d}=\pm1$;
\item $|S_2|=|S_1|+1$, the splitting circle from $D_{S_1}$ and both the
resulting circles from $D_{S_2}$ are marked with $\varb$, then $(S_1:S_2)_{\GF
d}=\pm1$.
\end{enumerate}
\end{enumerate}

Denote by $H(D)$ the homology with respect to $\GF+d$. It can be shown that
$H(D)$ is invariant under the Reidemeister moves, so that we can safely write
$H(L)$ instead. The fact that $p\not=2$ is crucial here, as the proof involves
division by 2 (see~\cite{Lee}).

Theorem~4.4 from~\cite{Lee} that states
\begin{equation}
H(L)\cong\frac{\Ker\bigl(\GF:\CalH_{\Z_p}(L)\to\CalH_{\Z_p}(L)\bigr)}
{\Im\bigl(\GF:\CalH_{\Z_p}(L)\to\CalH_{\Z_p}(L)\bigr)}
\end{equation}
still holds true without changes. One needs to use the fact that $L$ is
$\Z_p$H-thin, since it is H-slim, here.

The only non-trivial generalization is proving an analogue of Theorem~4.2
of~\cite{Lee} that $\dim_{\Z_p}H(L)=2^m$. Lee's proof uses Hodge
theory arguments which are not applicable to the $\Z_p$ case. Fortunately for
us, Hodge theory is only used to provide a {\em lower} bound on $\dim_\Q
H(L;\Q)$. It follows that $2^m\le\dim_\Q H(L;\Q)\le\dim_{\Z_p}H(L)$, where the
former inequality is provided by Theorem~4.2 from~\cite{Lee}, and
the latter one by the Universal Coefficient Theorem. In particular, all the
enhanced states of $D$, such that at every crossing the two touching circles
have different labels, are linearly independent in $H(L)$. Such states are in
one-to-one correspondence with all the orientations of $L$
(see~\cite{Lee}). Lee's proof of the fact that
$\dim_\Q H(L;\Q)\le2^m$ still works without changes for $\Z_p$. Hence
$2^m\le\dim_{\Z_p}H(L)\le2^m$ and $\dim_{\Z_p}H(L)=2^m$.

Filling the remaining technical gaps is left to the reader.
\end{proof}

\subsection{Proof of Theorem~\ref{thm:H-slim==>no-odd-torsion}}
Let $L$ be an H-slim link. Since
$h_{\Z_p}^{i,j}(L)=h^{i,j}(L)+T_p^{i,j}(L)+T_p^{i+1,j}(L)$, it follows from
Theorems~\ref{thm:lee-link-thin} and~\ref{thm:Zp-link-thin} that
$h_{\Z_p}(L)-h(L)$ are arranged in ``knight move'' pairs everywhere without
having to subtract anything (cf., for example,
Conjecture~\ref{conj:barnatan-knight}). Hence,
\begin{equation}
T_p^{i,2i-\Gs(L)-1}(L)=T_p^{i+1,2i-\Gs(L)+3}(L)+ T_p^{i+2,2i-\Gs(L)+3}(L).
\end{equation}
Since the support of the Khovanov homology is finite, all $T_p^{i,j}(L)$ must
be zero. \qed

\section*{Appendix}\def\thesection{A}\refstepcounter{section}

This section contains information about standard and reduced Khovanov homology
of knots whose torsion has some remarkable properties. The knot pictures below
were generated using Robert Scharein's program KnotPlot~\cite{KnotPlot}.

\subsection{How to read the tables}
Columns and rows of the tables below are marked with $i$- and $j$-grading of
the Khovanov homology, respectively. For the standard homology, the
$j$-grading is always odd and the corresponding table entries are printed in
boldface. The reduced homology have their $j$-grading even. They occupy places
between the main rows.

Only entries representing non-trivial groups are shown. An entry of the form
$a,b_2,c_4$ means that the corresponding group is
$\Z^a\oplus\Z_2^b\oplus\Z_4^c$.

\newpage
\subsection{The knot $8_{19}$}\label{app:8-19-knot}
This is the first $\Q$H-thick knot. It is T-fancy as well.
\begin{table}[h]
\centerline{%
\vbox{\hbox{\includegraphics[scale=0.25]{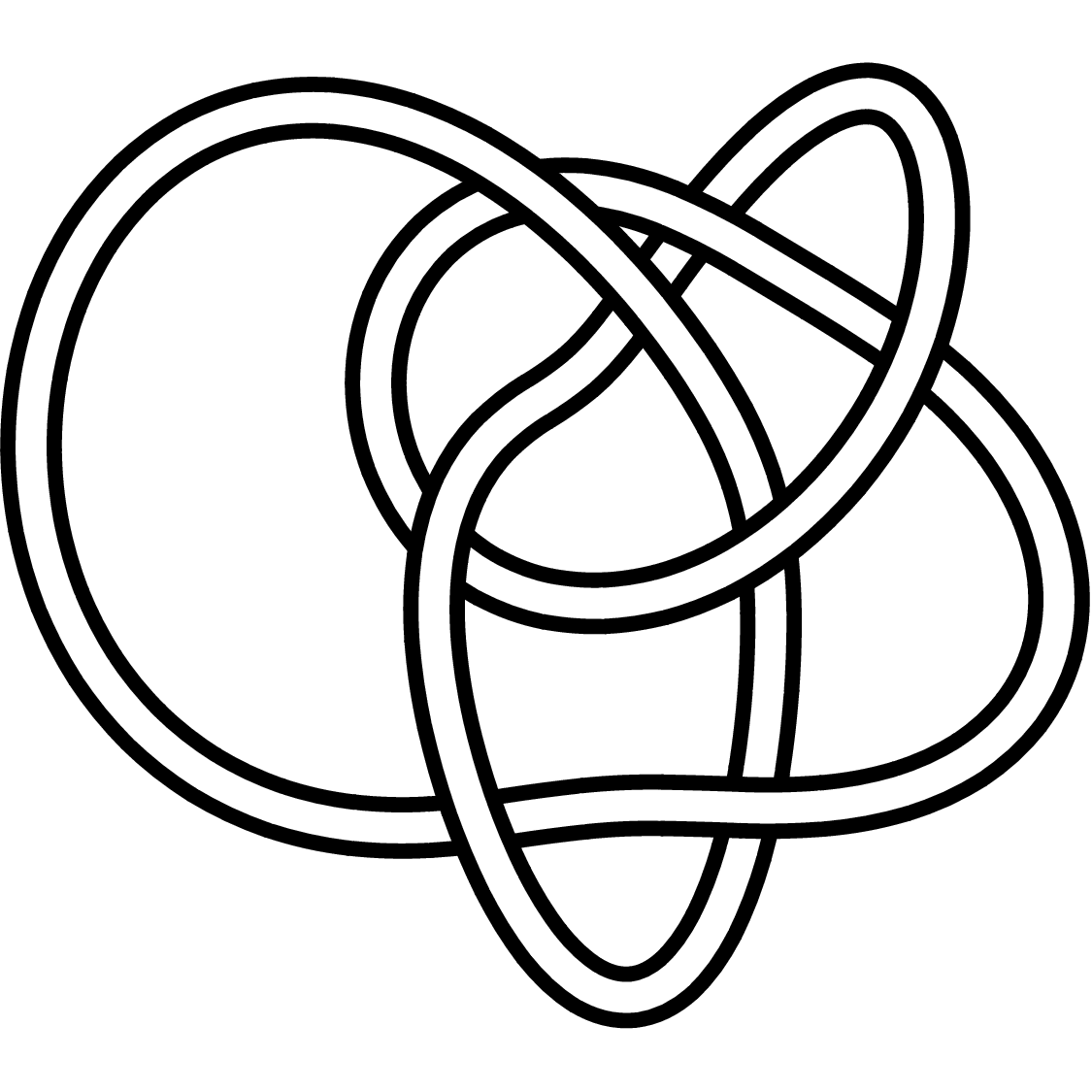}}\bigskip\null}%
\qquad\def\emptyline{&&&&&&&}
\def\numcolumns{9}%

\setbox\tablebox\vbox{\offinterlineskip\ialign{%
\vrule\TSp\vrule #\strut&\DSp\hfil #\DSp\vrule\TSp\vrule&
\DSp\hfil\bf #\hfil\DSp\vrule&
\DSp\hfil\bf #\hfil\DSp\vrule&
\DSp\hfil\bf #\hfil\DSp\vrule&
\DSp\hfil\bf #\hfil\DSp\vrule&
\DSp\hfil\bf #\hfil\DSp\vrule&
\DSp\hfil\bf #\hfil\DSp\vrule&#\TSp\vrule\cr
\dblhline
height 11pt depth 4pt&&
\rm\DSp0\DSp&
\rm\DSp1\DSp&
\rm\DSp2\DSp&
\rm\DSp3\DSp&
\rm\DSp4\DSp&
\rm\DSp5\DSp&\cr\dblhline
height 15pt depth 9pt&17&
&
&
&
&
&
\DSp1\DSp&
\cr
\omit\span\omit\leaders\hrule\hfill&\omit\leaders\hrule\hfill&
\omit\leaders\hrule\hfill&
\omit\leaders\hrule\hfill&
\omit\leaders\hrule\hfill&
\omit\leaders\hrule\hfill&
\omit\inline{1}&
\omit\leaders\hrule\hfill\vrule\cr
height 15pt depth 9pt&15&
&
&
&
&
&
\DSp1\DSp&
\cr
\omit\span\omit\leaders\hrule\hfill&\omit\leaders\hrule\hfill&
\omit\leaders\hrule\hfill&
\omit\leaders\hrule\hfill&
\omit\leaders\hrule\hfill&
\omit\leaders\hrule\hfill&
\omit\leaders\hrule\hfill&
\omit\leaders\hrule\hfill\vrule\cr
height 15pt depth 9pt&13&
&
&
&
\DSp1\DSp&
\DSp1\DSp&
&
\cr
\omit\span\omit\leaders\hrule\hfill&\omit\leaders\hrule\hfill&
\omit\leaders\hrule\hfill&
\omit\leaders\hrule\hfill&
\omit\inline{1}&
\omit\inline{1}&
\omit\leaders\hrule\hfill&
\omit\leaders\hrule\hfill\vrule\cr
height 15pt depth 9pt&11&
&
&
&
$\fam6{1}_{2}$&
\DSp1\DSp&
&
\cr
\omit\span\omit\leaders\hrule\hfill&\omit\leaders\hrule\hfill&
\omit\leaders\hrule\hfill&
\omit\inline{1}&
\omit\leaders\hrule\hfill&
\omit\leaders\hrule\hfill&
\omit\leaders\hrule\hfill&
\omit\leaders\hrule\hfill\vrule\cr
height 15pt depth 9pt&9&
&
&
\DSp1\DSp&
&
&
&
\cr
\omit\span\omit\leaders\hrule\hfill&\omit\leaders\hrule\hfill&
\omit\leaders\hrule\hfill&
\omit\leaders\hrule\hfill&
\omit\leaders\hrule\hfill&
\omit\leaders\hrule\hfill&
\omit\leaders\hrule\hfill&
\omit\leaders\hrule\hfill\vrule\cr
height 15pt depth 9pt&7&
\DSp1\DSp&
&
&
&
&
&
\cr
\omit\span\omit\leaders\hrule\hfill&\omit\inline{1}&
\omit\leaders\hrule\hfill&
\omit\leaders\hrule\hfill&
\omit\leaders\hrule\hfill&
\omit\leaders\hrule\hfill&
\omit\leaders\hrule\hfill&
\omit\leaders\hrule\hfill\vrule\cr
height 15pt depth 9pt&5&
\DSp1\DSp&
&
&
&
&
&
\cr
\dblhline
}}

\box\tablebox

}%
\medskip
\caption{The knot $8_{19}$ and its Khovanov homology (standard and reduced).}
\label{table:8-19-knot}
\end{table}

\subsection{The knot $9_{42}$}\label{app:9-42-knot}
This knot is $\Q$H-thick but T-simple. Conjecture~\ref{conj:H-T-simple} states
that every T-fancy knot should be H-thick as well. This example shows that the
converse is not true in general.
\begin{table}[h]
\centerline{%
\vbox{\hbox{\includegraphics[scale=0.3]{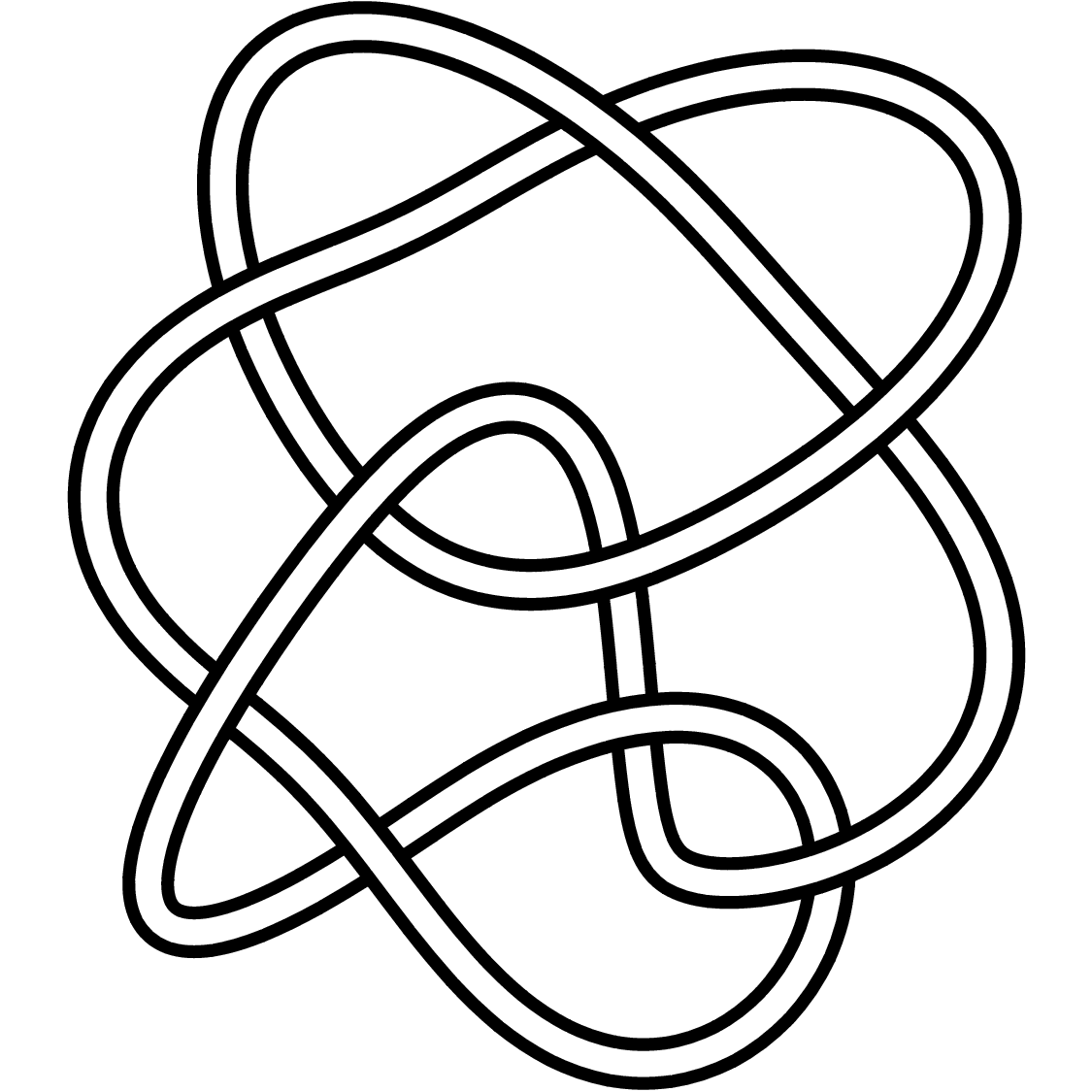}}\bigskip\null}%
\def\emptyline{&&&&&&&&}
\def\numcolumns{10}%

\setbox\tablebox\vbox{\offinterlineskip\ialign{%
\vrule\TSp\vrule #\strut&\DSp\hfil #\DSp\vrule\TSp\vrule&
\DSp\hfil\bf #\hfil\DSp\vrule&
\DSp\hfil\bf #\hfil\DSp\vrule&
\DSp\hfil\bf #\hfil\DSp\vrule&
\DSp\hfil\bf #\hfil\DSp\vrule&
\DSp\hfil\bf #\hfil\DSp\vrule&
\DSp\hfil\bf #\hfil\DSp\vrule&
\DSp\hfil\bf #\hfil\DSp\vrule&#\TSp\vrule\cr
\dblhline
height 11pt depth 4pt&&
\rm\DSp-4\DSp&
\rm\DSp-3\DSp&
\rm\DSp-2\DSp&
\rm\DSp-1\DSp&
\rm\DSp0\DSp&
\rm\DSp1\DSp&
\rm\DSp2\DSp&\cr\dblhline
height 15pt depth 9pt&7&
&
&
&
&
&
&
\DSp1\DSp&
\cr
\omit\span\omit\leaders\hrule\hfill&\omit\leaders\hrule\hfill&
\omit\leaders\hrule\hfill&
\omit\leaders\hrule\hfill&
\omit\leaders\hrule\hfill&
\omit\leaders\hrule\hfill&
\omit\leaders\hrule\hfill&
\omit\inline{1}&
\omit\leaders\hrule\hfill\vrule\cr
height 15pt depth 9pt&5&
&
&
&
&
&
&
$\fam6{1}_{2}$&
\cr
\omit\span\omit\leaders\hrule\hfill&\omit\leaders\hrule\hfill&
\omit\leaders\hrule\hfill&
\omit\leaders\hrule\hfill&
\omit\leaders\hrule\hfill&
\omit\leaders\hrule\hfill&
\omit\inline{1}&
\omit\leaders\hrule\hfill&
\omit\leaders\hrule\hfill\vrule\cr
height 15pt depth 9pt&3&
&
&
&
&
\DSp1\DSp&
\DSp1\DSp&
&
\cr
\omit\span\omit\leaders\hrule\hfill&\omit\leaders\hrule\hfill&
\omit\leaders\hrule\hfill&
\omit\leaders\hrule\hfill&
\omit\leaders\hrule\hfill&
\omit\inline{1}&
\omit\leaders\hrule\hfill&
\omit\leaders\hrule\hfill&
\omit\leaders\hrule\hfill\vrule\cr
height 15pt depth 9pt&1&
&
&
&
\DSp1\DSp&
1, $\fam6{1}_{2}$&
&
&
\cr
\omit\span\omit\leaders\hrule\hfill&\omit\leaders\hrule\hfill&
\omit\leaders\hrule\hfill&
\omit\leaders\hrule\hfill&
\omit\inline{2}&
\omit\inline{1}&
\omit\leaders\hrule\hfill&
\omit\leaders\hrule\hfill&
\omit\leaders\hrule\hfill\vrule\cr
height 15pt depth 9pt&-1&
&
&
&
1, $\fam6{1}_{2}$&
\DSp1\DSp&
&
&
\cr
\omit\span\omit\leaders\hrule\hfill&\omit\leaders\hrule\hfill&
\omit\leaders\hrule\hfill&
\omit\inline{1}&
\omit\leaders\hrule\hfill&
\omit\leaders\hrule\hfill&
\omit\leaders\hrule\hfill&
\omit\leaders\hrule\hfill&
\omit\leaders\hrule\hfill\vrule\cr
height 15pt depth 9pt&-3&
&
\DSp1\DSp&
\DSp1\DSp&
&
&
&
&
\cr
\omit\span\omit\leaders\hrule\hfill&\omit\leaders\hrule\hfill&
\omit\inline{1}&
\omit\leaders\hrule\hfill&
\omit\leaders\hrule\hfill&
\omit\leaders\hrule\hfill&
\omit\leaders\hrule\hfill&
\omit\leaders\hrule\hfill&
\omit\leaders\hrule\hfill\vrule\cr
height 15pt depth 9pt&-5&
&
$\fam6{1}_{2}$&
&
&
&
&
&
\cr
\omit\span\omit\leaders\hrule\hfill&\omit\inline{1}&
\omit\leaders\hrule\hfill&
\omit\leaders\hrule\hfill&
\omit\leaders\hrule\hfill&
\omit\leaders\hrule\hfill&
\omit\leaders\hrule\hfill&
\omit\leaders\hrule\hfill&
\omit\leaders\hrule\hfill\vrule\cr
height 15pt depth 9pt&-7&
\DSp1\DSp&
&
&
&
&
&
&
\cr
\dblhline
}}

\box\tablebox

}%
\medskip
\caption{The knot $9_{42}$ and its Khovanov homology (standard and reduced).}
\label{table:9-42-knot}
\end{table}

\subsection{The knot $13^n_{3663}$}\label{app:13n-3663-knot}
This is the first T-rich knot (the groups
with excessive torsion are $\CalH^{-3,-7}$, $\CalH^{-3,-5}$, $\CalH^{-2, -5}$,
$\CalH^{-2,-3}$, $\CalH^{0,-1}$, $\CalH^{0,1}$, $\CalH^{1,1}$, and
$\CalH^{1,3}$). This knot has $2$-torsion in the reduced homology as well.
This supports the claim of Conjecture~\ref{conj:reduced=T-rich} that a knot is
T-rich if and only if its reduced Khovanov homology has torsion. This knot is also the
first one whose homology is supported on 4 diagonals. The only other knots
with 13 crossings or less that share the same properties are $13^n_{4587}$,
$13^n_{4639}$, and $13^n_{5016}$.
\begin{table}[h]
\centerline{\def\emptyline{&&&&&&&&&&&&&&&}
\def\numcolumns{17}%

\setbox\tablebox\vbox{\offinterlineskip\ialign{%
\vrule\TSp\vrule #\strut&\DSp\hfil #\DSp\vrule\TSp\vrule&
\DSp\hfil\bf #\hfil\DSp\vrule&
\DSp\hfil\bf #\hfil\DSp\vrule&
\DSp\hfil\bf #\hfil\DSp\vrule&
\DSp\hfil\bf #\hfil\DSp\vrule&
\DSp\hfil\bf #\hfil\DSp\vrule&
\DSp\hfil\bf #\hfil\DSp\vrule&
\DSp\hfil\bf #\hfil\DSp\vrule&
\DSp\hfil\bf #\hfil\DSp\vrule&
\DSp\hfil\bf #\hfil\DSp\vrule&
\DSp\hfil\bf #\hfil\DSp\vrule&
\DSp\hfil\bf #\hfil\DSp\vrule&
\DSp\hfil\bf #\hfil\DSp\vrule&
\DSp\hfil\bf #\hfil\DSp\vrule&
\DSp\hfil\bf #\hfil\DSp\vrule&#\TSp\vrule\cr
\dblhline
height 11pt depth 4pt&&
\rm\DSp-6\DSp&
\rm\DSp-5\DSp&
\rm\DSp-4\DSp&
\rm\DSp-3\DSp&
\rm\DSp-2\DSp&
\rm\DSp-1\DSp&
\rm\DSp0\DSp&
\rm\DSp1\DSp&
\rm\DSp2\DSp&
\rm\DSp3\DSp&
\rm\DSp4\DSp&
\rm\DSp5\DSp&
\rm\DSp6\DSp&
\rm\DSp7\DSp&\cr\dblhline
height 15pt depth 9pt&13&
&
&
&
&
&
&
&
&
&
&
&
&
&
\DSp1\DSp&
\cr
\omit\span\omit\leaders\hrule\hfill&\omit\leaders\hrule\hfill&
\omit\leaders\hrule\hfill&
\omit\leaders\hrule\hfill&
\omit\leaders\hrule\hfill&
\omit\leaders\hrule\hfill&
\omit\leaders\hrule\hfill&
\omit\leaders\hrule\hfill&
\omit\leaders\hrule\hfill&
\omit\leaders\hrule\hfill&
\omit\leaders\hrule\hfill&
\omit\leaders\hrule\hfill&
\omit\leaders\hrule\hfill&
\omit\leaders\hrule\hfill&
\omit\inline{1}&
\omit\leaders\hrule\hfill\vrule\cr
height 15pt depth 9pt&11&
&
&
&
&
&
&
&
&
&
&
&
&
&
$\fam6{1}_{2}$&
\cr
\omit\span\omit\leaders\hrule\hfill&\omit\leaders\hrule\hfill&
\omit\leaders\hrule\hfill&
\omit\leaders\hrule\hfill&
\omit\leaders\hrule\hfill&
\omit\leaders\hrule\hfill&
\omit\leaders\hrule\hfill&
\omit\leaders\hrule\hfill&
\omit\leaders\hrule\hfill&
\omit\leaders\hrule\hfill&
\omit\leaders\hrule\hfill&
\omit\leaders\hrule\hfill&
\omit\leaders\hrule\hfill&
\omit\inline{1}&
\omit\leaders\hrule\hfill&
\omit\leaders\hrule\hfill\vrule\cr
height 15pt depth 9pt&9&
&
&
&
&
&
&
&
&
&
&
&
\DSp1\DSp&
\DSp1\DSp&
&
\cr
\omit\span\omit\leaders\hrule\hfill&\omit\leaders\hrule\hfill&
\omit\leaders\hrule\hfill&
\omit\leaders\hrule\hfill&
\omit\leaders\hrule\hfill&
\omit\leaders\hrule\hfill&
\omit\leaders\hrule\hfill&
\omit\leaders\hrule\hfill&
\omit\leaders\hrule\hfill&
\omit\leaders\hrule\hfill&
\omit\leaders\hrule\hfill&
\omit\leaders\hrule\hfill&
\omit\inline{1}&
\omit\leaders\hrule\hfill&
\omit\leaders\hrule\hfill&
\omit\leaders\hrule\hfill\vrule\cr
height 15pt depth 9pt&7&
&
&
&
&
&
&
&
&
&
\DSp1\DSp&
\DSp1\DSp&
$\fam6{1}_{2}$&
&
&
\cr
\omit\span\omit\leaders\hrule\hfill&\omit\leaders\hrule\hfill&
\omit\leaders\hrule\hfill&
\omit\leaders\hrule\hfill&
\omit\leaders\hrule\hfill&
\omit\leaders\hrule\hfill&
\omit\leaders\hrule\hfill&
\omit\leaders\hrule\hfill&
\omit\leaders\hrule\hfill&
\omit\leaders\hrule\hfill&
\omit\inline{1}&
\omit\inline{2}&
\omit\leaders\hrule\hfill&
\omit\leaders\hrule\hfill&
\omit\leaders\hrule\hfill&
\omit\leaders\hrule\hfill\vrule\cr
height 15pt depth 9pt&5&
&
&
&
&
&
&
&
&
\DSp1\DSp&
$\fam6{1}_{2}$&
1, $\fam6{1}_{2}$&
&
&
&
\cr
\omit\span\omit\leaders\hrule\hfill&\omit\leaders\hrule\hfill&
\omit\leaders\hrule\hfill&
\omit\leaders\hrule\hfill&
\omit\leaders\hrule\hfill&
\omit\leaders\hrule\hfill&
\omit\leaders\hrule\hfill&
\omit\leaders\hrule\hfill&
\omit\leaders\hrule\hfill&
\omit\inline{2}&
\omit\inline{1}&
\omit\leaders\hrule\hfill&
\omit\leaders\hrule\hfill&
\omit\leaders\hrule\hfill&
\omit\leaders\hrule\hfill&
\omit\leaders\hrule\hfill\vrule\cr
height 15pt depth 9pt&3&
&
&
&
&
&
&
&
$\fam6{1}_{2}$&
2, $\fam6{1}_{2}$&
\DSp1\DSp&
&
&
&
&
\cr
\omit\span\omit\leaders\hrule\hfill&\omit\leaders\hrule\hfill&
\omit\leaders\hrule\hfill&
\omit\leaders\hrule\hfill&
\omit\leaders\hrule\hfill&
\omit\leaders\hrule\hfill&
\omit\leaders\hrule\hfill&
\omit\leaders\hrule\hfill&
\omit\inline{1, ${1}_{2}$}&
\omit\inline{1}&
\omit\leaders\hrule\hfill&
\omit\leaders\hrule\hfill&
\omit\leaders\hrule\hfill&
\omit\leaders\hrule\hfill&
\omit\leaders\hrule\hfill&
\omit\leaders\hrule\hfill\vrule\cr
height 15pt depth 9pt&1&
&
&
&
&
&
&
2, $\fam6{1}_{2}$&
1, $\fam6{1}_{2}$&
$\fam6{1}_{2}$&
&
&
&
&
&
\cr
\omit\span\omit\leaders\hrule\hfill&\omit\leaders\hrule\hfill&
\omit\leaders\hrule\hfill&
\omit\leaders\hrule\hfill&
\omit\leaders\hrule\hfill&
\omit\leaders\hrule\hfill&
\omit\leaders\hrule\hfill&
\omit\inline{2, ${1}_{2}$}&
\omit\inline{1}&
\omit\leaders\hrule\hfill&
\omit\leaders\hrule\hfill&
\omit\leaders\hrule\hfill&
\omit\leaders\hrule\hfill&
\omit\leaders\hrule\hfill&
\omit\leaders\hrule\hfill&
\omit\leaders\hrule\hfill\vrule\cr
height 15pt depth 9pt&-1&
&
&
&
&
\DSp1\DSp&
\DSp1\DSp&
1, $\fam6{2}_{2}$&
\DSp1\DSp&
&
&
&
&
&
&
\cr
\omit\span\omit\leaders\hrule\hfill&\omit\leaders\hrule\hfill&
\omit\leaders\hrule\hfill&
\omit\leaders\hrule\hfill&
\omit\leaders\hrule\hfill&
\omit\inline{1}&
\omit\inline{2}&
\omit\leaders\hrule\hfill&
\omit\leaders\hrule\hfill&
\omit\leaders\hrule\hfill&
\omit\leaders\hrule\hfill&
\omit\leaders\hrule\hfill&
\omit\leaders\hrule\hfill&
\omit\leaders\hrule\hfill&
\omit\leaders\hrule\hfill&
\omit\leaders\hrule\hfill\vrule\cr
height 15pt depth 9pt&-3&
&
&
&
&
$\fam6{2}_{2}$&
1, $\fam6{1}_{2}$&
&
&
&
&
&
&
&
&
\cr
\omit\span\omit\leaders\hrule\hfill&\omit\leaders\hrule\hfill&
\omit\leaders\hrule\hfill&
\omit\leaders\hrule\hfill&
\omit\inline{1}&
\omit\inline{1, ${1}_{2}$}&
\omit\leaders\hrule\hfill&
\omit\leaders\hrule\hfill&
\omit\leaders\hrule\hfill&
\omit\leaders\hrule\hfill&
\omit\leaders\hrule\hfill&
\omit\leaders\hrule\hfill&
\omit\leaders\hrule\hfill&
\omit\leaders\hrule\hfill&
\omit\leaders\hrule\hfill&
\omit\leaders\hrule\hfill\vrule\cr
height 15pt depth 9pt&-5&
&
&
&
1, $\fam6{1}_{2}$&
1, $\fam6{1}_{2}$&
&
&
&
&
&
&
&
&
&
\cr
\omit\span\omit\leaders\hrule\hfill&\omit\leaders\hrule\hfill&
\omit\leaders\hrule\hfill&
\omit\leaders\hrule\hfill&
\omit\inline{${1}_{2}$}&
\omit\leaders\hrule\hfill&
\omit\leaders\hrule\hfill&
\omit\leaders\hrule\hfill&
\omit\leaders\hrule\hfill&
\omit\leaders\hrule\hfill&
\omit\leaders\hrule\hfill&
\omit\leaders\hrule\hfill&
\omit\leaders\hrule\hfill&
\omit\leaders\hrule\hfill&
\omit\leaders\hrule\hfill&
\omit\leaders\hrule\hfill\vrule\cr
height 15pt depth 9pt&-7&
&
\DSp1\DSp&
&
$\fam6{1}_{2}$&
&
&
&
&
&
&
&
&
&
&
\cr
\omit\span\omit\leaders\hrule\hfill&\omit\leaders\hrule\hfill&
\omit\inline{1}&
\omit\leaders\hrule\hfill&
\omit\leaders\hrule\hfill&
\omit\leaders\hrule\hfill&
\omit\leaders\hrule\hfill&
\omit\leaders\hrule\hfill&
\omit\leaders\hrule\hfill&
\omit\leaders\hrule\hfill&
\omit\leaders\hrule\hfill&
\omit\leaders\hrule\hfill&
\omit\leaders\hrule\hfill&
\omit\leaders\hrule\hfill&
\omit\leaders\hrule\hfill&
\omit\leaders\hrule\hfill\vrule\cr
height 15pt depth 9pt&-9&
&
$\fam6{1}_{2}$&
&
&
&
&
&
&
&
&
&
&
&
&
\cr
\omit\span\omit\leaders\hrule\hfill&\omit\inline{1}&
\omit\leaders\hrule\hfill&
\omit\leaders\hrule\hfill&
\omit\leaders\hrule\hfill&
\omit\leaders\hrule\hfill&
\omit\leaders\hrule\hfill&
\omit\leaders\hrule\hfill&
\omit\leaders\hrule\hfill&
\omit\leaders\hrule\hfill&
\omit\leaders\hrule\hfill&
\omit\leaders\hrule\hfill&
\omit\leaders\hrule\hfill&
\omit\leaders\hrule\hfill&
\omit\leaders\hrule\hfill&
\omit\leaders\hrule\hfill\vrule\cr
height 15pt depth 9pt&-11&
\DSp1\DSp&
&
&
&
&
&
&
&
&
&
&
&
&
&
\cr
\dblhline
}}

\box\tablebox
}
\medskip
\caption{Standard and reduced Khovanov homology of the knot $13^n_{3663}$.}
\label{table:13n-3663-knot}
\end{table}

\newpage
\subsection{The $(4,5)$-torus knot}\label{app:4-5-torus}
This is one of the first knots whose Khovanov homology has torsion of
order~$4$. Its minimal diagram has 15 crossings. There are no knots with 13
crossings or less that have torsion of order other than~$2$. This knot is also
T-rich and has $2$-torsion in reduced homology (cf.
Conjecture~\ref{conj:reduced=T-rich}).
\begin{table}[h]
\centerline{\def\emptyline{&&&&&&&&&&&&}
\def\numcolumns{14}%

\setbox\tablebox\vbox{\offinterlineskip\ialign{%
\vrule\TSp\vrule #\strut&\DSp\hfil #\DSp\vrule\TSp\vrule&
\DSp\hfil\bf #\hfil\DSp\vrule&
\DSp\hfil\bf #\hfil\DSp\vrule&
\DSp\hfil\bf #\hfil\DSp\vrule&
\DSp\hfil\bf #\hfil\DSp\vrule&
\DSp\hfil\bf #\hfil\DSp\vrule&
\DSp\hfil\bf #\hfil\DSp\vrule&
\DSp\hfil\bf #\hfil\DSp\vrule&
\DSp\hfil\bf #\hfil\DSp\vrule&
\DSp\hfil\bf #\hfil\DSp\vrule&
\DSp\hfil\bf #\hfil\DSp\vrule&
\DSp\hfil\bf #\hfil\DSp\vrule&#\TSp\vrule\cr
\dblhline
height 11pt depth 4pt&&
\rm\DSp0\DSp&
\rm\DSp1\DSp&
\rm\DSp2\DSp&
\rm\DSp3\DSp&
\rm\DSp4\DSp&
\rm\DSp5\DSp&
\rm\DSp6\DSp&
\rm\DSp7\DSp&
\rm\DSp8\DSp&
\rm\DSp9\DSp&
\rm\DSp10\DSp&\cr\dblhline
height 15pt depth 9pt&29&
&
&
&
&
&
&
&
&
&
&
$\fam6{1}_{2}$&
\cr
\omit\span\omit\leaders\hrule\hfill&\omit\leaders\hrule\hfill&
\omit\leaders\hrule\hfill&
\omit\leaders\hrule\hfill&
\omit\leaders\hrule\hfill&
\omit\leaders\hrule\hfill&
\omit\leaders\hrule\hfill&
\omit\leaders\hrule\hfill&
\omit\leaders\hrule\hfill&
\omit\leaders\hrule\hfill&
\omit\leaders\hrule\hfill&
\omit\inline{${1}_{2}$}&
\omit\leaders\hrule\hfill\vrule\cr
height 15pt depth 9pt&27&
&
&
&
&
&
&
&
&
&
\DSp1\DSp&
$\fam6{1}_{2}$&
\cr
\omit\span\omit\leaders\hrule\hfill&\omit\leaders\hrule\hfill&
\omit\leaders\hrule\hfill&
\omit\leaders\hrule\hfill&
\omit\leaders\hrule\hfill&
\omit\leaders\hrule\hfill&
\omit\leaders\hrule\hfill&
\omit\leaders\hrule\hfill&
\omit\leaders\hrule\hfill&
\omit\leaders\hrule\hfill&
\omit\inline{1}&
\omit\leaders\hrule\hfill&
\omit\leaders\hrule\hfill\vrule\cr
height 15pt depth 9pt&25&
&
&
&
&
&
&
&
\DSp1\DSp&
&
$\fam6\torframe{$\fam6{1}_{4}$}$&
&
\cr
\omit\span\omit\leaders\hrule\hfill&\omit\leaders\hrule\hfill&
\omit\leaders\hrule\hfill&
\omit\leaders\hrule\hfill&
\omit\leaders\hrule\hfill&
\omit\leaders\hrule\hfill&
\omit\leaders\hrule\hfill&
\omit\leaders\hrule\hfill&
\omit\inline{1}&
\omit\inline{1}&
\omit\leaders\hrule\hfill&
\omit\leaders\hrule\hfill&
\omit\leaders\hrule\hfill\vrule\cr
height 15pt depth 9pt&23&
&
&
&
&
&
\DSp1\DSp&
&
1, $\fam6{1}_{2}$&
\DSp1\DSp&
&
&
\cr
\omit\span\omit\leaders\hrule\hfill&\omit\leaders\hrule\hfill&
\omit\leaders\hrule\hfill&
\omit\leaders\hrule\hfill&
\omit\leaders\hrule\hfill&
\omit\leaders\hrule\hfill&
\omit\inline{1}&
\omit\leaders\hrule\hfill&
\omit\inline{${1}_{2}$}&
\omit\leaders\hrule\hfill&
\omit\leaders\hrule\hfill&
\omit\leaders\hrule\hfill&
\omit\leaders\hrule\hfill\vrule\cr
height 15pt depth 9pt&21&
&
&
&
&
&
\DSp1\DSp&
\DSp1\DSp&
$\fam6{1}_{2}$&
&
&
&
\cr
\omit\span\omit\leaders\hrule\hfill&\omit\leaders\hrule\hfill&
\omit\leaders\hrule\hfill&
\omit\leaders\hrule\hfill&
\omit\leaders\hrule\hfill&
\omit\leaders\hrule\hfill&
\omit\leaders\hrule\hfill&
\omit\inline{1}&
\omit\leaders\hrule\hfill&
\omit\leaders\hrule\hfill&
\omit\leaders\hrule\hfill&
\omit\leaders\hrule\hfill&
\omit\leaders\hrule\hfill\vrule\cr
height 15pt depth 9pt&19&
&
&
&
\DSp1\DSp&
\DSp1\DSp&
&
\DSp1\DSp&
&
&
&
&
\cr
\omit\span\omit\leaders\hrule\hfill&\omit\leaders\hrule\hfill&
\omit\leaders\hrule\hfill&
\omit\leaders\hrule\hfill&
\omit\inline{1}&
\omit\inline{1}&
\omit\leaders\hrule\hfill&
\omit\leaders\hrule\hfill&
\omit\leaders\hrule\hfill&
\omit\leaders\hrule\hfill&
\omit\leaders\hrule\hfill&
\omit\leaders\hrule\hfill&
\omit\leaders\hrule\hfill\vrule\cr
height 15pt depth 9pt&17&
&
&
&
$\fam6{1}_{2}$&
\DSp1\DSp&
&
&
&
&
&
&
\cr
\omit\span\omit\leaders\hrule\hfill&\omit\leaders\hrule\hfill&
\omit\leaders\hrule\hfill&
\omit\inline{1}&
\omit\leaders\hrule\hfill&
\omit\leaders\hrule\hfill&
\omit\leaders\hrule\hfill&
\omit\leaders\hrule\hfill&
\omit\leaders\hrule\hfill&
\omit\leaders\hrule\hfill&
\omit\leaders\hrule\hfill&
\omit\leaders\hrule\hfill&
\omit\leaders\hrule\hfill\vrule\cr
height 15pt depth 9pt&15&
&
&
\DSp1\DSp&
&
&
&
&
&
&
&
&
\cr
\omit\span\omit\leaders\hrule\hfill&\omit\leaders\hrule\hfill&
\omit\leaders\hrule\hfill&
\omit\leaders\hrule\hfill&
\omit\leaders\hrule\hfill&
\omit\leaders\hrule\hfill&
\omit\leaders\hrule\hfill&
\omit\leaders\hrule\hfill&
\omit\leaders\hrule\hfill&
\omit\leaders\hrule\hfill&
\omit\leaders\hrule\hfill&
\omit\leaders\hrule\hfill&
\omit\leaders\hrule\hfill\vrule\cr
height 15pt depth 9pt&13&
\DSp1\DSp&
&
&
&
&
&
&
&
&
&
&
\cr
\omit\span\omit\leaders\hrule\hfill&\omit\inline{1}&
\omit\leaders\hrule\hfill&
\omit\leaders\hrule\hfill&
\omit\leaders\hrule\hfill&
\omit\leaders\hrule\hfill&
\omit\leaders\hrule\hfill&
\omit\leaders\hrule\hfill&
\omit\leaders\hrule\hfill&
\omit\leaders\hrule\hfill&
\omit\leaders\hrule\hfill&
\omit\leaders\hrule\hfill&
\omit\leaders\hrule\hfill\vrule\cr
height 15pt depth 9pt&11&
\DSp1\DSp&
&
&
&
&
&
&
&
&
&
&
\cr
\dblhline
}}

\box\tablebox
}
\medskip
\caption{Standard and reduced Khovanov homology of the $(4,5)$-torus knot.}
\label{table:4-5-torus}
\end{table}

\raggedright

\end{document}